\documentclass[11pt]{amsart}                          
\usepackage{epsfig}       
\usepackage{subfig}     
\usepackage{amssymb,amscd,bbm}     
\usepackage{amsthm,amsgen,amsmath,epic,psfrag}
\usepackage{rotating}
\usepackage[dvips]{xypic}       

\xyoption{curve}                
\xyoption{rotate}               

\newcommand{\pa}{{\mathcal P}}

\newcommand{\R}{{\mathbb R}}

\newcommand{\F}{{\mathbb F}}

\newcommand{\Z}{\mathbb{Z}}

\newcommand{\conf}{{\rm Conf}}
\newcommand{\uconf}{\overline{\rm Conf}}
\newcommand{\Sym}{{\si}}
\newcommand{\sm}{\wedge}

\newcommand{\K}{{\mathcal{K}}}
\newcommand{\mD}{{\mathcal{D}}}
\newcommand{\orb}{{\rm Orb}}
\def\cupp{\cdot}
\newcommand{\tr}{\odot}

\newcommand{\und}{\underline}

\newcommand{\ka}{{\bf k}}
\newcommand{\sym}{{\rm Sym}}
\newcommand{\scr}{{\mathcal Q}}
\newcommand{\si}{{\mathcal S}}
\newcommand{\A}{\mathcal A}
\newcommand{\sq}{{Sq}}
\newcommand{\sqi}{{Sq^{i}}}
\newcommand{\rep}{{\rm Rep}}

\theoremstyle{plain}
\newtheorem{theorem}{Theorem}[section]
\newtheorem{proposition}[theorem]{Proposition}
\newtheorem{lemma}[theorem]{Lemma}
\newtheorem{corollary}[theorem]{Corollary}

\theoremstyle{definition}
\newtheorem{definition}[theorem]{Definition}

\newtheorem{example}[theorem]{Example}

\theoremstyle{remark}
\newtheorem{remark}[theorem]{Remark}

\newcommand{\refT}[1]{Theorem~\ref{T:#1}}
\newcommand{\refC}[1]{Corollary~\ref{C:#1}}

\newcommand{\refP}[1]{Proposition~\ref{P:#1}}
\newcommand{\refD}[1]{Definition~\ref{D:#1}}
\newcommand{\refL}[1]{Lemma~\ref{L:#1}}
\newcommand{\refE}[1]{Equation~\ref{E:#1}}

\input{colordvi}

\begin{document}

\title{The mod-two cohomology rings of symmetric groups}
\author[C. Giusti]{Chad Giusti}
\address{Mathematics Department, 
Willamette University}
\email{cgiusti@willamette.edu}
\author[P. Salvatore]{Paolo Salvatore}
\address{Dipartimento di Matematica, � 
Universit\`a di Roma �Tor Vergata�}
\email{salvator@mat.uniroma2.it}
\author[D. Sinha]{Dev Sinha}
\address{Mathematics Department, 
University of Oregon}
\email{dps@math.uoregon.edu}

\begin{abstract}
We present a new additive basis for the mod-two cohomology of symmetric groups, along with explicit
rules for multiplication and application of Steenrod operations in that basis.  The key organizational tool is
a Hopf ring structure introduced by Strickland and Turner.  We elucidate some of the relationships between 
our approach and previous approaches to the homology and cohomology of symmetric groups.
\end{abstract}

\maketitle

\section{Introduction}

We determine the mod-two 
cohomology of all symmetric groups, that is of the disjoint union 
$\coprod_{n}  B\si_{n}$,
 as a Hopf ring.  This description allows us to give the first 
 additive basis with a complete, explicit rule for multiplication in that basis.
 We also give geometric representatives for mod-two cohomology and
 explicitly describe the  action of the Steenrod algebra .

\begin{definition}\label{D:Hopfring}
A Hopf ring is  a ring object in the category of coalgebras.  Explicitly, a Hopf ring is
vector space $V$ with two multiplications, one 
comultiplication, and an antipode $(\tr, \cdot, \Delta, S)$ such that the first multiplication forms a Hopf algebra
with the comultiplication and antipode, the second multiplication forms a bialgebra with the comultiplication,
and these structures satisfy the distributivity relation 
\begin{equation*}
\alpha \cdot (\beta \tr \gamma) =  \sum_{\Delta \alpha = \sum a' \otimes a''} (a' \cdot \beta) \tr 
(a'' \cdot \gamma).
\end{equation*}
\end{definition}

We consider only Hopf rings where all of these structures are commutative.

On the cohomology of $\coprod_{n} B \si_{n}$ the second product $\cdot$ is cup product, 
which is zero for classes supported on disjoint components.
The first product $\tr$ is the  transfer product  - see  \refD{odot} -
first studied by Strickland and Turner \cite{StTu97}.  
It is akin to the ``induction product'' in the representation theory of symmetric
groups  \cite{Knut73, Zele81}.
The coproduct $\Delta$ on cohomology is dual to the standard Pontrjagin product
on the homology of $\coprod_{n} B \si_{n}$.   

\begin{theorem}\label{T:genandrel}
As a Hopf ring, $H^{*}(\coprod_{n} B \si_{n}; \F_{2})$ is generated by 
classes $\gamma_{\ell, n} \in H^{n(2^{\ell} - 1)}(B \si_{n 2^{\ell}})$,
along with unit classes on each component.
The coproduct of $\gamma_{\ell, n}$ is given by
$$\Delta \gamma_{\ell, n} = \sum_{i+j = n} {\gamma_{\ell, i}} \otimes {\gamma_{\ell, j}}.$$
Relations between transfer products of these generators are given by 
$$\gamma_{\ell, n} \tr \gamma_{\ell, m} = \binom{n+m}{n} \gamma_{\ell, n+m}.$$
The antipode is the identity map.
Cup products of generators
on different components are zero, and there are no other relations between cup products of generators.
\end{theorem}

Thus, all of the relations in the cohomology of symmetric groups follow from the distributivity of
cup product over transfer product. 
Building on this   presentation we  give  an additive basis,
which is fairly immediate, and an explicit presentation of multiplication rules 
in that basis.  This additive basis is represented graphically by ``skyline diagrams'' which
are reminiscent of Young diagrams.  The rule for cup product is  complicated but
accessible, akin to rules for multiplying symmetrized monomials.

We begin the paper by developing Hopf rings which arise in algebra related to those
we study, namely that of symmetric
invariants and in representations of symmetric groups.
Though Hopf rings were introduced by Milgram to study the homology of the sphere
spectrum \cite{Milg70} and thus the infinite symmetric group \cite{BaPr72},  and later used to 
study other ring spectra \cite{RaWi76},
the Hopf ring structure we study does not fit into that framework.  In particular
it exists in cohomology rather than homology.  See \cite{StTu97} for a lucid
discussion of the relationships between all of these structures.

We show that  these Hopf ring generators are, and thus all cohomology is, represented by 
Thom classes of linear subvarieties.
We connect with previous work and identify the restriction maps in cohomology
to elementary abelian subgroups.  
We  use such restriction maps to study the action of the Steenrod algebra.  There is a Cartan 
formula for the transfer product, so the Steenrod
action on the cohomology of symmetric groups is completely determined by that on the 
Hopf ring generators $\gamma_{\ell,2^{k}}$, which we give in \refT{steen}.

We  revisit some of Feshbach's calculations \cite{Fesh02} and
express his cup-product generators in terms of our Hopf ring generators.  
While the 
Hopf ring presentation of all components is  straightforward, the cup ring structure for a single symmetric group is still complicated.  
We also give our own invariant-theoretic presentation.  
At the end of the paper we show that  
 Stiefel-Whitney classes for the standard representations can be used as Hopf ring generators, 
 forging another tie between the categories of finite sets and vector spaces.

The cohomology of symmetric groups is a classical topic, dating back to Steenrod's \cite{Stee53} and Adem's \cite{Adem57}
studies of them in the context of cohomology operations.  We heavily rely on 
Nakaoka's seminal work \cite{Naka61} which determined the mod-two homology of symmetric groups.  
More explicit treatment of the cup product structure on cohomology 
was later given at the prime two partially by Hu'ng \cite{Hung87} and Adem-McGannis-Milgram \cite{AMM90}
and more definitively by
Feshbach \cite{Fesh02}, using
restriction to elementary abelian subgroups and invariant theory.
 While Feshbach's generators for cup ring structure are accessible, the relations
are given recursively, in increasingly complex forms.
The Hopf ring structure give a  compact, closed-form recursive description of all
components at once.  It seems that it will also be useful other primes, to other groups, to other
configuration spaces, and to related spaces.
 
\tableofcontents
We thank Nick Kuhn for pointing out a simpler proof of Theorem~\ref{T:primitive} 
as well as Nick Proudfoot, Hal Sadofsky and Alejandro Adem for helpful conversations.  
  The third author would 
like to thank the Universities of Roma Tor Vergata, Pisa, Zurich, Lille, Louvain and Nice, 
and the CIRM  in France and the INDAM agency in Italy for their hospitality.

\section{Hopf rings arising from representations and invariants of symmetric groups}

In this section we identify some Hopf rings defined by classical objects, namely  
rings of invariants  and representation rings 
of symmetric groups, which are related to the
Hopf ring structure on the cohomology of symmetric groups.  
Though invariant theory and representation theory
have long, distinguished histories, to our knowledge
the use of Hopf rings to serve as a framework for restriction and induction maps is new.

\begin{definition}\label{D:invtHopf}
Let $A$ be an algebra which is flat over a ground ring $R$ (which is suppressed from
notation).  Let $\mu_{m,n} : A^{\otimes m} \otimes A^{\otimes n} \to A^{\otimes m+n}$
denote the standard isomorphism, and let $\Delta_{m,n}$ denote its inverse.

Let $A^{\si} = \bigoplus_{n} (A^{\otimes n})^{\si_{n}}$, which 
we call the total symmetric invariants of $A$.  Define a 
coproduct $\Delta$ to be the sum of restrictions of $\Delta_{m,n}$.  Define a product
$\tr : (A^{\otimes m})^{\si_{m}} \otimes (A^{\otimes n})^{\si_{n}} \to
(A^{\otimes m+n})^{\si_{m+n}}$ as the symmetrization  of $\mu_{m,n}$ 
over $\si_{m+n}/ (\si_{m} \times \si_{n})$.

If $A$ is a polynomial algebra, define an antipode $S$ on $A^{\otimes n}$ which multiplies a monomial
by  $(-1)^{k}$, where $k$ the number of variables which appear
in the monomial.
\end{definition}

For explicit calculations with symmetric invariants, which we make throughout this section, we set the following
notation.

\begin{definition}
Let $\sym(m)$ denote the minimal symmetrization of a monomial $m$, namely $\sum_{[\sigma] \in S_{n}/H} \;\;  \sigma \cdot m$ 
where $H$ is the subgroup which fixes $m$.  
\end{definition}

\begin{proposition}\label{P:invtHopf}
The total symmetric invariants of $A$, namely $A^{\si}$, with the product $\tr$,
its standard product (which is zero for elements from
different summands),  and the coproduct $\Delta$ forms a Hopf semiring.
When $A$ is a polynomial algebra, the total symmetric invariants forms a Hopf ring
with the antipode $S$.
\end{proposition}

\begin{remark}
This construction can be generalized in significant ways.  First, the 
rings $A^{\otimes n}$ can be replaced by more general rings with 
$\si_{n}$ action and analogues of maps $\mu$ and $\Delta$.
More generally, they could be replaced by schemes, obtaining Hopf
rings through regular functions or perhaps  some sort of cohomology.
Also, instead of symmetric groups other sequences of groups with
inclusions $G_{n} \times G_{m} \to G_{m+n}$, in particular linear groups over
finite fields, can be used.  We content ourselves here with the minimum needed
to treat cohomology of symmetric groups.
\end{remark}

\begin{proof}[Proof of \refP{invtHopf}]
The fact that the standard product and $\Delta$ form a bialgebra follows from the fact that the 
$\Delta_{m,n}$ are ring homomorphisms.  

That $\tr$ and $\Delta$ form a bialgebra is also possible to establish for all of 
$A^{\otimes m} \otimes A^{\otimes n}$, not
just the $\si_{m} \times \si_{n}$-invariants.  First we consider 
$a_{1} \otimes \cdots \otimes a_{m}  \tr b_{1} \otimes \cdots \otimes b_{n} $,
which by definition the symmetrization over $\si_{m+n}$ of 
$\tau = a_{1} \otimes \cdots \otimes a_{n}  \otimes b_{1} \otimes \cdots \otimes b_{m}$.  We choose
this symmetrization to be given by shuffles.  
Next we  apply $\Delta$ and consider the  $\Delta_{m',n'}$-summand, which
takes the first  $m'$ and last $n'$ tensor factors of a given tensor.
The result of applying $\Delta_{m',n'}$ to one of the shuffles at hand will be a
 shuffle of $a_{1}, \ldots, a_{i}$ and $b_{1}, \ldots, b_{j}$
tensored with a shuffle of $a_{i+1}, \ldots, a_{m}$ with $b_{j+1}, \ldots, b_{n}$.  But these
pairs of smaller shuffles 
are exactly what is obtained if one first applies $\Delta_{i,m-i} \otimes \Delta_{j,n-j}$ to $\tau$
and then $\tr$-multiplies, establishing the result.

For distributivity we start with $a \in (A^{\otimes m+n})^{\si_{m+n}}$ and $b$ and $c$
in $A^{\otimes m}$ and $A^{\otimes n}$ respectively.  Then $a \cdot (b \tr c)$ is the product
of $a$ with the symmetrization by shuffles of $b \otimes c$.  But since $a$ is already symmetric 
this is equal to the  symmetrization of $a \cdot (\mu_{m,n} b\otimes c) =
\mu_{m,n} (\Delta_{m,n}(a) \cdot b \otimes c)$.  Since there are no other terms in the 
coproduct of $a$ which non-trivially multiply $b \otimes c$, we get that 
$a \cdot (b \tr c) = \sum_{\Delta a = a' \otimes a''} a' \cdot b \tr a'' \tr c$.

 Now we restrict to when  $A$ is a polynomial algebra generated by some $\{ x_{i} \}$, 
in which case the antipode map $S$ multiplies a monomial by 
 $(-1)^{k}$, where $k$ the number of variables which appear in the monomial.
Consider symmetrizations of the form
$$\sym({x_{P}}^{\vec{q}}) = 
\sym\left( {x_1}^{q_{1}}{x_{{2}}}^{q_{1}}  \cdots {x_{p_{1}}}^{q_{1}} {x_{p_{1} + 1}}^{q_{2}} \cdots {x_{p_{1} + p_{2}}}^{q_{2}}
\cdots {x_{n - p_{k}}}^{q_{k}} \cdots {x_{n}}^{q_{k}}\right),$$
where $P$ is the partition  $\sum p_{i} = n$.  These 
span the symmetric invariants of $A$ so we check that $S$ is an antipode by applying
$\mu_{\tr}  \circ (S \otimes id) \circ \Delta$.
The coproduct of $\sym({x_{P}}^{\vec{q}})$ is $\sum_{P = P' + P''} \sym({x_{P'}}^{\vec{q}}) \otimes \sym({x_{P''}}^{\vec{q}})$,
where $P'$ and $P''$ vary over partitions with $p_{i}  = p_{i}' + p_{i}''.$  Applying $S \times id$ introduces a sign of 
$(-1)^{|P'|}$ to each term in the sum, where $|P'| = \sum p_{i}'.$ 
Applying $\tr$, each term $ \sym({x_{P'}}^{\vec{q}}) \tr \sym({x_{P''}}^{\vec{q}})$ produces a multiple of 
of the original symmetrized monomial $\sym({x_{P}}^{\vec{q}})$.  
Because $\sym({x_{P}}^{\vec{q}})$ had $\frac{n!}{P!}$ terms 
(where  for a partition $P!$ is the product of $p_{i}!$),  and 
$ \sym({x_{P'}}^{\vec{q}}) \tr \sym({x_{P''}}^{\vec{q}})$ has $\frac{n!}{P'! P''!}$ terms, this multiple is  
$(-1)^{|P'|} \frac{P!}{P'! P''!}  \sym({x_{P}}^{\vec{q}})$.  That $S$ is an antipode then follows from the identify 
$\sum_{P = P' + P''}   (-1)^{|P'|} \frac{P!}{P'! P''!}  = 0,$ which generalizes the familiar fact for binomial coefficients when
$P$ is a singleton partition.

In summary, what we have proven is that both $(\cdot, \Delta)$ and $(\tr, \Delta)$ define
bialgebra structures on all of $\oplus_{n} A^{\otimes n}$, which then restrict to invariants.  
Moreover, distributivity
of $\cdot$ over $\tr$ holds when multiplying something which is $\si_{n}$
invariant, which means that when  restricting to the total symmetric invariants we obtain a Hopf semiring.
Finally, when $A$ is a polynomial algebra $(\tr, \Delta, S)$ define a Hopf algebra structure on the total
symmetric invariants, so we obtain a Hopf ring.
\end{proof}

\begin{definition}\label{D:Hopfmono}
A Hopf ring monomial in classes $x_{i}$ is one of the form $f_{1} \tr f_{2} \tr \cdots \tr f_{k}$, where
each $f_{j}$ is a monomial under the $\cdot$ product in the $x_{i}$.
\end{definition}

These monomials play a significant role in all of our examples.

\begin{example}\label{Ex:classical}
The total symmetric invariants $\ka[x]^{\si}$ is, as a vector space,  
the direct sum of the classical rings of symmetric 
polynomials over $\ka$.
We do not know whether this Hopf ring structure on the direct sum of all symmetric
polynomials has been considered previously.  

The second product in the Hopf ring structure is the standard product of symmetric polynomials,
defined to be zero if the number of variables differs.
The coproduct is  ``de-coupling'' of two sets of variables followed by reindexing, so for example
$$\Delta_{2,1}({x_{1}}^{2} x_{2} x_{3} + {x_{1}} {x_{2}}^{2} x_{3} + {x_{1}} x_{2} {x_{3}}^{2})
= ({x_{1}}^{2} x_{2} + x_{1}{x_{2}}^{2}) \otimes x_{1} + x_{1}x_{2}\otimes {x_{1}}^{2}.$$

The first product $f \tr g$ reindexes the variables of $g$, multiplies that by $f$, and then
symmetrizes with respect to $\si_{n+m}/ \si_{n} \times \si_{m},$ as
can be done with shuffles.  A reasonable name for this product would be the
shuffle product.  For example
$$({x_{1}}^{2}x_{2} + x_{1}{x_{2}}^{2}) \tr x_{1} =
({x_{1}}^{2}x_{2} + x_{1}{x_{2}}^{2})x_{3} + ({x_{1}}^{2}x_{3} + x_{1}{x_{3}}^{2})x_{2}
+ ({x_{2}}^{2}x_{3} + x_{2}{x_{3}}^{2})x_{1} = 2 \sym({x_{1}}^{2} x_{2} x_{3}).$$

Let $1_{k}$ denote the unit function on $k$ variables
and $\sigma_{n}(k)$ the $n$th symmetric function in $k$ variables.  Because $\sigma_{n}(k)
= \sigma_{n}(n) \tr 1_{k-n}$, 
as a Hopf ring symmetric functions are generated by the $1_{k}$ and 
$\sigma_{n} = \sigma_{n}(n)$.
The $\sigma_{n}$ $\tr$-multiply according to the rule $\sigma_{n} \tr \sigma_{m}
= \binom{n+m}{n} \sigma_{n+m}$, a divided powers algebra.

Thus over the rationals only $\sigma_{1}$ is required to generate as a Hopf ring, 
while over $\F_{p}$ one needs all $\sigma_{p^{i}}$.    Because of periodicity
of binomial coefficients modulo $p$, the Hopf sub-rings generated by classes
$\sigma_{n p^{i}}$ for fixed $i$ (or equivalently, the quotients  obtained by 
setting other symmetric polynomials to zero) are isomorphic 
to the full Hopf ring of symmetric functions.  This isomorphism accounts for some
``self-similarity'' in the cohomology of symmetric groups.

Irreducible Hopf ring monomials in the $\sigma_{n}$ correspond to symmetrized monomials in 
the $x_{i}$.  That is, for $p_{i}$ distinct,
\begin{equation}\label{E:sym}
{\sigma_{n_{1}}}^{p_{1}} \tr {\sigma_{n_{2}}}^{p_{2}} \tr \cdots \tr {\sigma_{n_{k}}}^{p_{k}} \tr 1_{j}
= \sym \left( {x_{1}}^{p_{1}} \cdots {x_{n_{1}}}^{p_{1}} {x_{n_{1}+ 1}}^{p_{2}} \cdots {x_{n_{1}+ n_{2}}}^{p_{2}}
\cdots {x_{N}}^{p_{k}} \right),
\end{equation}
where $N =\sum n_{i}$.
Distributivity in the Hopf ring structure gives rise to an inductive method
to multiply symmetrized monomials.  
This approach to symmetric polynomials is 
fairly indifferent to the classical theorem that the ring of symmetric polynomials in a fixed number of 
variables forms a polynomial algebra.
\end{example}

More generally, we may let 
$A =  \ka[x_{\bf 1}, \ldots, x_{\bf m}]$ in which case 
total symmetric invariants  are the direct sum of rings of symmetric polynomials in $m$ collections of variables.
Explicitly, we take polynomials in variables $x_{{\bf i}, j}$ with $1 \leq i \leq m$ and $1 \leq j \leq n$ which
are invariant under permutation of the second subscripts (so that the bold subscripts are ``fixed'').
Define $\sigma_{{\bf i}, n}$ to be  $x_{{\bf i}, 1} \cdot x_{{\bf i}, 2}  \cdot \cdots \cdot x_{{\bf i}, n}$. 

\begin{example}\label{Ex:multivar}
Consider  $\ka[x_{{\bf 1}, 1}, x_{{\bf 1}, 2}, x_{{\bf 2}, 1}, x_{{\bf 2}, 2}]^{\si_{2}}$.  Key elements are
\begin{enumerate}
\item $\sigma_{{\bf 1},1} \tr 1 = x_{{\bf 1}, 1} + x_{{\bf 1}, 2}$
\item $\sigma_{{\bf 1},2}  = x_{{\bf 1}, 1} x_{{\bf 1}, 2}$
\item $\sigma_{{\bf 2},1} \tr 1 = x_{{\bf 2}, 1} + x_{{\bf 2}, 2}$
\item $\sigma_{{\bf 2},2}  = x_{{\bf 2}, 1} x_{{\bf 2}, 2}$
\item $\sigma_{{\bf 1}, 1} \sigma_{{\bf 2}, 1} \tr 1 = x_{{\bf 1}, 1} x_{{\bf 2},1} + x_{{\bf 1}, 2} x_{{\bf 2},2}$.
\item $\sigma_{{\bf 1}, 1} \tr \sigma_{{\bf 2}, 1}  = x_{{\bf 1}, 1} x_{{\bf 2},2} + x_{{\bf 1}, 2} x_{{\bf 2},1}$.
\end{enumerate}

\end{example}

We then have the following.

\begin{proposition}\label{P:symmfun}
The total symmetric invariants of  $A =  \ka[x_{{\bf i}, j}]$ is generated as
a Hopf ring by unit elements and the elementary products $\sigma_{{\bf i}, n}$.  
The coproduct is given by 
$$\Delta \sigma_{{\bf i},n} = \sum_{j+k = n} \sigma_{{\bf i}, j} \otimes \sigma_{{\bf i}, k}.$$
The $\tr$-products are given by $$\sigma_{{\bf i},  n} \tr \sigma_{{\bf i}, m} = 
\binom{n+m}{n} \sigma_{{\bf i}, n+m},$$ while $\tr$-products between classes with
different $\ell$ are free.    The collection of $\sigma_{{\bf i},n}$  for all $i$ with fixed $n$ form a polynomial ring
under the standard product.
\end{proposition}

In this presentation, the standard product is determined
by the relations given above, 
Hopf ring distributivity, the fact that products of classes in different rings of invariants are zero,
and the fact that the collection of $\sigma_{{\bf i},n}$ with fixed $n$ form a polynomial ring.  
As we discuss in Section~\ref{S:invt}, when $\ka = {\mathbb F}_{2}$ these Hopf rings are isomorphic to split 
quotient-Hopf rings of the cohomology of symmetric groups.

\begin{proof}
That these symmetric functions have coproducts, $\tr$-products and ordinary products as
stated is straightforward.  The fact that these basic symmetric functions in each
collection of variables are Hopf ring generators follows from the fact that their associated 
Hopf monomial basis coincides with the symmetrized monomial basis for the symmetric
polynomials.  The case of one set of variables is given in Equation~\ref{E:sym} above.
More generally we translate between these two bases as follows.

To translate from the Hopf monomial basis to symmetric polynomials proceeds as already defined.
The simple monomial ${\sigma_{{\bf i_{1}}, n}}^{p_{1}} \cdot \cdots \cdot {\sigma_{{\bf i_{k}}, n}}^{p_{k}}$ 
(note that the number of variables $n$ for each symmetric function must be the same to have the product non-zero) 
is by definition 
$$\left( {x_{{\bf i_{1}}, 1}}^{p_{1}} \cdots {x_{{\bf i_{1}}, n}}^{p_{1}} \right) \cdot \cdots \cdot 
\left( { x_{{\bf i_{k}}, 1}}^{p_{k}} \cdots {x_{{\bf i_{k}}, n}}^{p_{k}} \right).$$
Hopf monomials are transfer products of these, which by definition translate the second indices and then symmetrize.
Let us denote by $\phi$ this map of sets from the set of Hopf ring monomials  to the set of symmetrized monomials in $\ka[x_{{\bf i}, j}]$.
For example, $\phi( \sigma_{{\bf 3}, 2} \tr {\sigma_{{\bf 1},1}}^{3} {\sigma_{{\bf 3}, 1}}^{3} \tr {\sigma_{{\bf 1}, 5}}^{7}) $ 
is by definition equal to 
$$\sym\left( (x_{{\bf 3},1} x_{{\bf 3},2}) ({x_{{\bf 1}, 3}}^{3} {x_{{\bf 3}, 3}}^{3}) ( {x_{{\bf 1}, 4}}^{7} \cdot \cdots \cdot {x_{{\bf 1}, 8}}^{7}) \right).$$

Conversely, we may start with the minimal symmetrization of an arbitrary monomial in the $x_{{\bf i}, j}$, which if we collect terms which share the same second index is of the form  
$$\sym \left( {x_{{\bf i_{1}}, j_{1}}}^{p_{1}} \cdot {x_{{\bf i_{2}}, j_{1}}}^{p_{2}} \cdots {x_{{\bf i_{m}}, j_{k}}}^{p_{m}} \cdots   
{x_{{\bf i_{\ell}}, j_{k}}}^{p_{\ell}} \right).$$
Call the factor of this product which shares a second index $j_{q}$ a $j$-factor.  We say two $j$-factors  are similar if they 
differ only by relabeling of this second index.  Then for each $j$-factor
${x_{{\bf i_{m}, j}}}^{p_{m}} \cdots {x_{{\bf i_{n}, j}}}^{p_{n}}$ and all of the $j$-factors which are similar to it we associate
the product ${\sigma_{{\bf i_{m}}, s}}^{p_{m}} \cdots {\sigma_{{\bf i_{n}}, s}}^{p_{n}},$ where $s$ is the number of 
$j$ factors which are similar.  We then form a Hopf monomial by taking $\tr$-products of these, over the similarity classes
of $j$-factors.  Let us denote by $\psi$ this map of sets from the set of symmetrized monomials in $\ka[x_{{\bf i}, j}]$  to the set of Hopf ring monomials of the total symmetric invariants of 
$\ka[x_{{\bf i}, j}]$.

For example, consider $z = \sym\left( {x_{{\bf 4}, 6}}^{23} {x_{{\bf 2}, 7}}^{18} {x_{{\bf 3}, 7}}^{5} {x_{{\bf 2}, 8}}^{12} {x_{{\bf 4}, 9}}^{23}\right).$
The $j$-factors ${x_{{\bf 4}, 6}}^{23}$ and ${x_{{\bf 4}, 9}}^{23}$ are similar, so to them we associate 
${\sigma_{{\bf 4}, 2}}^{23}$.  Thus $\psi(z) =  {\sigma_{{\bf 4}, 2}}^{23} \tr {\sigma_{{\bf 2}, 1}}^{18} {\sigma_{{\bf 3}, 1}}^{5} \tr {\sigma_{{\bf 2}, 1}}^{12}.$

We claim that $\phi \circ \psi$ is the identity.  Indeed, $\phi(\psi(\omega))$ differs from $\omega$ by the action of an element of the symmetric group which sends all of the second indices in the first $j$-factor chosen to $1, \cdots, s$, all of the second indices
in the second $j$-factor to $s+1, \cdots$, and so forth.  Since symmetrized monomials span all symmetric polynomials, this shows that Hopf monomials in the $\sigma_{{\bf i}, n}$ span the total symmetric invariants.

\end{proof}

As we will see as well in the case of cohomology of symmetric groups, this simple presentation of the Hopf ring structure belies the fact that understanding the standard product structure alone is complicated.  

\begin{example}
Consider two sets of
two variables, namely  $\ka[x_{{\bf 1}, 1}, x_{{\bf 1}, 2}, x_{{\bf 2}, 1}, x_{{\bf 2}, 2}]^{\si_{2}}$ 
as in Example~\ref{Ex:multivar}.  
The ring of invariants has an additive basis of ${\sigma_{{\bf 1}, 1}}^{p}  {\sigma_{{\bf 2}, 1}}^{q} \tr {\sigma_{{\bf 1}, 1}}^{r}
{\sigma_{{\bf 2}, 1}}^{s}$, with either $p \neq q$ or $r \neq s$, along with ${\sigma_{{\bf 1}, 2}}^{p} {\sigma_{{\bf 2}, 2}}^{r}$, 
which we can multiply using the Hopf ring structure.

Understanding this ring in terms of generators and relations is already involved.   There is a fourth-degree relation, namely
\begin{multline*}
0 = (\sigma_{{\bf 1}, 1}\sigma_{{\bf 2}, 1} \tr 1)^{2} + (\sigma_{{\bf 1}, 1}\sigma_{{\bf 2}, 1} \tr 1)  (\sigma_{{\bf 1}, 1} \tr 1) 
 (\sigma_{{\bf 2}, 1} \tr 1) +   \sigma_{{\bf 1}, 2}  (\sigma_{{\bf 2}, 1} \tr 1)  +  (\sigma_{{\bf 1}, 1} \tr 1)  \sigma_{{\bf 2}, 2}. 
\end{multline*}
\end{example}

Indeed, the classical theorem that symmetric functions in one set of variables form a
polynomial algebra is an anomaly, as even the simplest cases of multiple sets of variables
are quite involved.  The structure of such rings over $\F_{2}$ is the computational heart
of Adem-McGannis-Milgram and Feshbach's work on symmetric groups \cite{AMM90, Fesh02}.
To our knowledge, even generators of such rings of invariants have not been
computed over $\F_{p}$ with $p$ odd.   

\bigskip

Though we do not apply them in this paper, we take a moment to develop Hopf ring
structures on  representation rings of symmetric groups, since they are
direct analogues of the one we study in cohomology.

In his book \cite{Zele81}, Zelevinsky shows that the direct sum 
$\bigoplus \rep(\si_{n})$
forms a bialgebra, under the induction product (which he credits to Young) and
restriction coproduct.
Denote the induction product, which takes $V$ a representation of $\si_{n}$ and
$W$ a representation of $\si_{m}$ to ${\rm Ind}_{\si_{n} \times \si_{m}}^{\si_{n+m}} V \otimes W$,  
by $\tr$.  Denote the coproduct, which is the sum of maps which sends $V$ to 
${\rm Res}_{\si_{n} \times \si_{m}}^{\si_{n+m}}(V)$, 
by $\Delta$, to be consistent with notation from topology.  Denote the standard product 
in the representation ring, given by tensor product of representations, by $\cdot$. 
Let $x_{n}$ denote the trivial representation of $\si_{n}$, and define an antipode 
$S$ on these by setting $(\sum x_{i}) \tr \left(\sum S(x_{i}) \right) = 0$.

\begin{proposition}\label{P:rephopfring}
Over any ground field, $\bigoplus_{n} \rep(\si_{n})$ with induction product, 
tensor product (defined to be zero if representations are
of different symmetric groups) and restriction coproduct, forms a Hopf ring. That is,
both induction/restriction and $S$ form a Hopf algebra,  tensor/restriction defines a bialgebra, 
and there is a distributivity condition 
$$V \cdot (W_{1} \tr W_{2}) = \sum_{\Delta V = \sum V_{1} \otimes V_{2}} (V_{1} \cdot W_{1}) \tr
(V_{2} \cdot W_{2}).$$
\end{proposition}

\begin{proof}[Sketch of proof]
The proof that one obtains bialgebra structures, 
after unraveling definitions, follows from basic theorems on induction and restriction.
Zelevinsky did not consider the antipode as part of his definition of Hopf algebra (and we 
do not know whether it has been considered before, or know of a more natural construction).
But Zelevinsky did prove that for complex representations
the bialgebra defined by induction product and restriction coproduct alone 
is isomorphic to the polynomial algebra $\Z[x_{1},  x_{2}, \ldots]$, with coproduct
that $\Delta x_{n} = \sum x_{i} \otimes x_{j}$.  
For a Hopf algebra which as an algebra is such a polynomial ring, the antipode $S$ is unique as given.
\end{proof}

The change-of-basis
between the monomial basis in the $x_{i}$, which correspond to permutation representations induced
up from the trivial representation of block subgroups, and the basis of irreducible representations is
highly non-trivial.

Our Hopf ring structure along with Zelevinsky's theorem gives a classically known method to compute
tensor products of representations induced up from the trivial representation of block subgroups, that is of
permutation representations.
The Hopf monomial basis in this case
coincides with the monomial basis under $\tr$, since each ${x_{i}}^{2} = x_{i}$.   If $P$ is a partition,
namely $p_{1} + \cdots + p_{k} = n$, we let $x_{P}$ denote $x_{p_{1}} \tr \cdots \tr x_{p_{k}}$, 
the permutation representation of $S_{n}$ induced up from the trivial representation of $S_{P}$.

For example, for $\si_{3}$ we have as an additive basis: $x_{3}$ which is the trivial representation,
$x_{1 + 2}$ which is the standard representation, and $x_{1 + 1 + 1}$ which is the
regular representation.  We can compute using the Hopf ring formalism for example that
\begin{equation*}
 x_{1+2} \otimes x_{1 + 2} = (x_{1} \tr x_{2}) \cdot (x_{1} \tr x_{2}) 
= {x_{1}}^{2} \tr {{x_2}}^{ 2} + {{x_1}}^{2} \tr \left(({x_1} \tr {x_1}) \cdot x_{2}
\right) = x_{1} \tr x_{2} + {x_{1}}\tr x_{1} \tr x_{1}.
\end{equation*}

As mentioned, the  distributivity formula encodes the induction-restriction formula, and thus leads to a classical method
of computing tensor products of these permutation representations.  If $P$ and $Q$ are
partitions of $n$ then consider any matrix $\hat{A}$ with nonnegative integer entries such that the entries
of $i$th row of $A$ add up to $p_{i}$ and those of the $j$th column of $A$ add up to $q_{j}$.  Then the entries of $\hat{A}$
form another partition of $n$, which we call $A$ and say that $A$ is a product-refinement of $P$ and $Q$.  For example
if $P = Q = 1 + 2$ then two possibilities for $\hat{A}$ are $\begin{pmatrix} 1 & 0 \\ 0 & 1 \end{pmatrix}$ and 
$\begin{pmatrix} 0 & 1 \\ 1 & 1 \end{pmatrix}$.

The following classical theorem (see Example I.7.23(e) of  \cite{Macd95})
is straightforward to establish using Hopf ring distributivity.  We state it because it is 
a direct analogue of our description of multiplication in the cohomology of 
symmetric groups through an additive basis, given in \refT{generalcup}.

\begin{proposition}
If $x_{P}$ and $x_{Q}$ are permutation representations then $x_{P} \otimes x_{Q} \cong \bigoplus_{A} x_{A},$
where the sum is over $A$ which are product-refinements of $P$ and $Q$.
\end{proposition}

\bigskip

More basically, we can consider the direct sum of Burnside rings of symmetric groups $\bigoplus
A(\si_{n})$.  As usual $A(G)$ is the ring obtained by group completing the monoid of 
$G$-sets, and thus is the  representation ring of $G$ in the category of finite sets.   
We define multiplication between different summands to be zero.  Define coproduct and induction
product analogously to how they were defined for representations in vector spaces.  Once again
we obtain a Hopf ring.
For any group we can map $A(G)$ to $\rep(G)$ by using a $G$-set as a basis for a vector space. 
Collecting these maps gives a map $\bigoplus_{n}A(\si_{n}) \to \bigoplus_{n}{\rep}(\si_{n})$ which respects Hopf ring structures.
For symmetric groups these maps are surjective.

\section{Definition of the transfer product}

The classifying space for symmetric groups is often modeled by unordered configuration
spaces, which are a natural context to define the second product in our Hopf ring structure.
Let $\conf_{n}(X) = \{ (x_{1}, \ldots, x_{n}) \in X^{\times n} | x_{i} \neq x_{j} \; {\rm if} \; i \neq j \}$.
Let $\uconf_{n}(X) = \conf_{n}(X) / \si_{n}$, where $\si_{n}$ acts on $\conf_{n}(X)$ by 
permuting indices.
\begin{definition}\label{D:odot}
Consider the following maps
$$
\begin{CD}
\uconf_{m,n}(X) @>p>> \uconf_{n}(X) \times \uconf_{m}(X) \\
@VfVV \\
\uconf_{m+n}(X).
\end{CD}
$$
Here $\uconf_{m,n}(X)$ is the space of $m+n$ distinct points in $X$, $m$ of which
have one color and  $n$ of which have another.  The map $f$ forgets labels, and is a 
covering map with $\binom{n+m}{n}$ sheets.  The map $p$ is the product of maps which project onto 
each of  the two groups of points separately.

Define the transfer product $\tr$ as the the composite $\tau_{f} \circ p^{*}$,
where $\tau_{f}$ denotes the transfer map associated to $f$ on cohomology.
\end{definition}

When $X = \R^{\infty}$ so that $\uconf_{n}(X) \simeq B\si_{n}$, 
this product was previously studied by Strickland and Turner \cite{StTu97}.
In this case the map  $p$ is a homotopy equivalence, so this composite
is essentially the transfer map itself.  Moreover, in this case the map 
$f$  is homotopic to the map defining the product on $\coprod_{n} B \si_{n}$.
Recall that by either applying the classifying space functor to the standard inclusions
$\si_{n} \times \si_{m} \hookrightarrow \si_{n+m}$ or by taking
unions of unordered configurations in the $\uconf_{n}(\R^{\infty})$ model, we get a product
on $\coprod_{n} B \si_{n}$ which passes to a commutative product $*$ on its homology.
Its dual $\Delta$ defines
a cocommutative coalgebra structure on cohomology.

The following is immediate from  Theorem~3.2 of \cite{StTu97}. 

\begin{theorem}
The transfer product $\tr$  along with the cup product $\cupp$
and the coproduct $\Delta$ define a Hopf semiring structure on
$H^{*}(\coprod_{n} B \si_{n})$ with coefficients in any ring.  With mod-two coefficients, the identity map gives an antipode which 
defines a full Hopf ring structure.
\end{theorem}

\begin{proof}[Sketch of proof.]
That cup product and the coproduct $\Delta$ form a bialgebra follows immediately from the fact that $\Delta$ is induced by 
the covering map $f$ of \refD{odot}.  The Hopf ring distributivity follows similarly from the fact that $\tr$ is induced by the transfer
associated to $f$.  

That $\tr$ and $\Delta$ form a bialgebra is essentially a double-coset formula.   As Adams notes in \cite{Adam78} 
such formulae usually follow from naturality of transfer maps.  Start with the following pull-back diagram of covering maps 
$$
\begin{CD}
\bigsqcup \uconf_{p,q,r,s}(\R^{\infty}) @>e>> \uconf_{m,n}(\R^{\infty}) \\
@VfVV @VgVV \\
\bigsqcup \uconf_{i,j}(\R^{\infty})  @>h>>   \uconf_{m+n}(\R^{\infty}),
\end{CD}
$$ 
Here $\uconf_{p,q,r,s}(\R^{\infty})$ is defined as configurations of colored points, $p$ of which have one color,
$q$ of which have a second color, etc., and the first union is over indices such that
$p+q = m$, $r+s = n$, $p + r = i$ and $q + s = j$, while the second union is over those with $i + j = m + n$.  
All of the covering maps $e$, $f$, $g$ and $h$
merge or forget colors.
If $\alpha \in H^{*}(B\si_{n})$ and $\beta \in H^{*}(B\si_{m})$ then by definition
 $\Delta (\alpha \tr \beta) = h^{*} \circ g^{!} (\alpha  \otimes \beta)$.  But because this is a pullback square of covering spaces,
 this is equal to $f^{!} \circ e^{*} (\alpha \otimes \beta)$, which is $\Delta(\alpha) \tr \Delta(\beta)$.

Finally, for the antipode we need to pass from the space-level divided powers construction $\bigvee_{n} B\si_{n} =
\bigvee_{n} ES_{n} \ltimes_{S_{n}} (S^{0})^{\wedge n}$, which is called $DS^{0}$, to the spectrum version.
(See \cite{BMMS86} for an explanation of  the extension of the functor $D$ to spectra).  
In \cite{StTu97} the authors define the antipode using the additive inverse map on the spectrum $S^{0}$.
But in mod-two homology and cohomology, this map induces the identity.
\end{proof}

The antipode for other coefficient systems, which we do not consider here, 
will be similar to that given for symmetric polynomials in \refP{invtHopf}, as one can see from the connections
we develop in Section~\ref{S:invt}.

\begin{remark} \label{stable}
The transfer product is induced by a stable map, namely the transfer 
$$\tau_f: \Sigma^\infty \uconf_{m+n}(X) \longrightarrow \Sigma^\infty \uconf_{m,n}(X).$$ The composite

$$\Psi=\Sigma^\infty p \circ \tau_f    
 : \Sigma^\infty(\uconf_{n}(X)) \to \Sigma^\infty(\left(\uconf_{k}(X) \times \uconf_{n-k}(X)\right) $$
is a stable map 
inducing on homology the coproduct $\Delta_\tr$ dual to the transfer product.

Because these structures are all defined by applying cohomology to maps
and stable maps, the generalized cohomology of symmetric
groups for any ring theory forms a Hopf ring, as  established in Theorem~3.2 of
\cite{StTu97}.  Strickland in  \cite{Stri98} uses this to study
the Morava $E$-theory of symmetric groups.  
The $K$-theory of symmetric groups is a 
completion of the representation ring, by the Atiyah-Segal 
theorem \cite{AtSe69}, and the Hopf ring structure defined by Strickland-Turner agrees with that 
of \refP{rephopfring}.
\end{remark}




The cohomology and representation theory
of various types of linear groups over finite fields also form Hopf rings, as do some rings
of invariants under these groups such as Dickson algebras, as do the 
cohomology of some symmetric products.  Using these Hopf ring structures for further study
is likely to be fruitful.

Because the fundamental geometry underlying the transfer product is that of taking a 
configuration and partitioning it into two configurations in all possible ways, we sometimes
call it the partition product.  Indeed, partitioning is part of the geometry as seen
through Poincar\'e duality.
The usual cup product 
corresponds to intersection of Poincar\'e duals (that is, supports of representing Thom classes), 
which means taking the the locus of configurations of 
$n$ points which satisfy the conditions defining both of 
the two cocycles in question.  The locus defining the transfer product is similar, but we instead require that some $k$ points satisfy the first condition and
then the complementary $n-k$ points satisfy the latter condition.



\section{Review of homology of symmetric groups}

We now focus on calculations with $\F_{2}$ coefficients, recollecting standard facts to set notation.

\begin{definition}
Let $\K$ be the associative algebra over $\F_{2}$, with product $\circ$, generated by $q_{0}, q_{1}, \ldots$  with the following relations,
$$
{\rm (Adem)}  \;\; {\rm For} \;\;  m > n, \;\; 
q_{m} \circ q_{n} = \sum_{i} \binom{i - n - 1}{2i - m - n} q_{m + 2n - 2i} \circ q_{i}.
$$
Given a sequence $I = i_{1}, \cdots, i_{k}$ of non-negative integers, let $q_{I} = q_{i_{1}} \circ \cdots
\circ q_{i_{k}}$.  Using the Adem relations, $\K$ is spanned by $q_{I}$ whose entries are  non-decreasing.  We call such an $I$ {\em admissible}.   If such an $I$ has no zeros we call it {\em strongly admissible}.
\end{definition}

Following \cite{BiJo97}, we call $\K$ the Kudo-Araki algebra, to distinguish it from a closely related presentation usually called the Dyer-Lashof algebra.  The algebra $\K$ is one of the main characters in algebraic topology because it acts on the homology of any infinite loop space, or more generally any $E_{\infty}$-space (see I.1 of \cite{CLM76}).

\begin{definition} \label{defq}
An action of $\K$ on a graded algebra $A$ with product denoted $*$ and grading denoted $\deg$ is a map from $\K \otimes A \to A$, typically written using operational notation, with the following properties:
\begin{itemize}
\item (Action) $(q_{i} \circ q_{j})(a) = q_{i}(q_{j}(a))$.
\item (Grading)  $\deg q_{i} (a) = 2\deg a + i$.
\item (Squaring) $q_{0} (a) = a^{*2}$.  
\item  (Vanish) $q_{i}(1) = 0$ for $i > 0$.
\item (Cartan)  For any $a, b$, we have  $q_{n}(a * b) = \sum_{i+j = n} q_{i} (a) * q_{j} ( b)$.
\end{itemize}
If $\K$ acts on $A$ we call $A$ a $\K$-algebra.
\end{definition}

We denote powers in such an algebra $A$ by $a* \cdots *a = a^{*n}$.

As is standard, there is a free $\K$-algebra functor, left adjoint to the forgetful functor from $\K$-algebras to vector spaces, with
which  we can give the simplest reformulation of Nakaoka's seminal result. 

\begin{theorem}[\cite{Naka61}]\label{T:Nakaoka}
$H_{*}(\coprod_{n} B \si_{n})$, with its standard product $*$, is isomorphic to the free $\K$-algebra generated by $H_{0}(B \si_{1})$.  

Thus, as a ring under $*$ it is isomorphic to the polynomial algebra generated by the nonzero class  $\iota \in H_{0}(B \si_{1})$ and $q_{I}(\iota) \in H_{*}(B \si_{2^{k}})$ for $I$ strongly admissible.
\end{theorem}

The second statement, which is essentially Nakaoka's formulation, follows straightforwardly from the first statement.  We will often abuse notation and refer to $q_{I}(\iota)$ as simply $q_{I}
\in H_{|I|}B\si_{2^{k}}$ with $|I| = i_{1} + 2 i_{2} + \cdots + 2^{k-1}i_{k}$.  

We describe geometric representatives of the classes $q_I$, and then some of their dual cohomology classes.

\begin{definition}\label{D:QI}
Given $I = i_{1}, i_{2}, \cdots, i_{k}$ inductively define manifolds ${\orb}_{I}$ and maps
$Q_{I} : {\orb}_{I} \to \uconf_{2^{k}}( \R^{d})$ where $d > i_{\ell}$ for all $\ell$ as follows.
\begin{itemize}
\item If $I$ is empty ${\orb}_{I}$ is a point.  Otherwise, $\orb_{I} = S^{i_{1}} \times_{\Z/2} (\orb_{I'} 
\times \orb_{I'})$, the quotient of $\Z/2$ acting antipodally on $S^{i_{1}}$ and by permuting the two factors of $\orb_{I'}$, where $I' = i_{2}, \cdots, i_{k}$.
\item Let $\varepsilon = \frac{1}{4}$.  If $I$ is empty, $Q_{I}$ sends $\orb_{I}$ to $0 = \R^{0}$.
Otherwise, $Q_{I}\left (v \times_{\Z/2} (o_{1}, o_{2})\right)$ is given by 
$(v + \varepsilon Q_{I'}(o_{1})) \bigcup (-v + \varepsilon Q_{I'}(o_{2}))$.  Here we consider $v \in S^{i_{1}}$ to be a unit vector in $\R^{i_{1}+ 1}$.  The configuration  $v + \varepsilon Q_{I'}(o_{1})$ is the configuration obtained by scaling each point in $Q_{I'}(o_{1})$ by $\varepsilon$ and then adding $v$, perhaps after either the configuration or $v$ is included (canonically) into the larger of the two Euclidean spaces in which they are defined. 
\end{itemize}
\end{definition}

\begin{proposition}
The class $q_{I}$ in \refT{Nakaoka} is equal to $(Q_{I})_{*}[\orb_{I}] \in
H_{*}(\uconf_{2^{k}}(\R^{\infty}))$, where $[\orb_{I}]$ is the fundamental class of $\orb_{I}$.
\end{proposition}

\begin{figure}
$$\includegraphics[width=4cm]{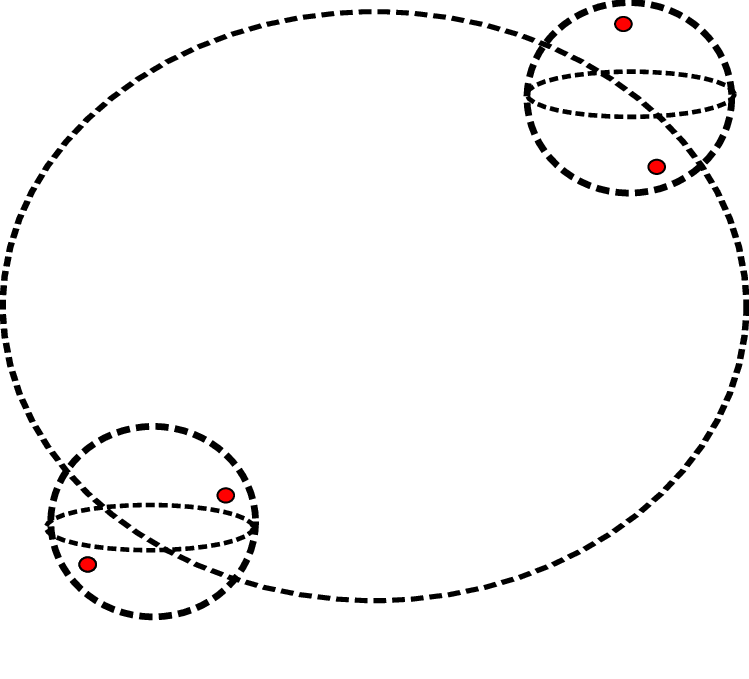}$$
\caption{An illustration of $q_{1,2} \in H_{5}(B\si_{4})$.}
\end{figure}

We now present ``linear'' geometric representatives for cohomology classes.  

\begin{definition}
Let $\hat{X}$ be a manifold without boundary of dimension $nd - m$ which maps properly to $\uconf_{n}(\R^{d})$.
We define its Thom class as follows.  Take the fundamental class in 
mod-two locally-finite homology of $\hat{X}$ in dimension $nd - m$, and map it to locally-finite homology
of $\uconf_{n}(\R^{d})$.  Apply the Poincar\'e duality isomorphism to obtain a class in
in $H^{m}(\uconf_{n}(\R^{d}); \F_{2})$, which we call the Thom class of $\hat{X}$.  By abuse we may sometimes refer
only to the image of $\hat{X}$.

If we choose $d$ large enough (greater than $m$), then
the restriction map from the cohomology of $\uconf_{n}(\R^{\infty}) = B\si_{n}$ to that of $\uconf_{n}(\R^{d})$ is an isomorphism
in degree $m$.  (One can deduce this  isomorphism 
from calculations in homology, which for $\uconf_{n}(\R^{d})$ is constructed from $q_{I}$ with $i_{k} < d$.)
So this Thom class lifts uniquely to define a class in the cohomology of $B\si_{n}$ which we also
call the Thom class.
\end{definition}

For example, if we refer to points
in $\uconf_{n}(\R^{d})$ as ${\bf x} = (x_{1}, \ldots, x_{n})/\sim$, then the non-zero class in 
$H^{1}(\uconf_{n}(\R^{d}))$ for $d \geq 2$
is represented by the variety $X$ of points such that some $x_{i}$ and $x_{j}$
share their first coordinate.  Specifically, let $\hat{X}$ be
the space of configurations of $n$ points, two of which have one color - say black - and the rest of 
which share another color, such that the two black points must share their first coordinate.
The Thom class of $\hat{X}$ mapping to the configuration space by forgetting colors is the non-trivial
class in degree one, as we can see by evaluating it on the cycle $q_{1}$ by intersection (exactly once
do two unlabeled points which are antipodal on some generic $S^{1}$ share their first coordinate).

Such Thom classes represent Hopf ring generators of the mod-two cohomology of symmetric groups.

 \begin{definition}\label{D:gammas}
Let $q_{\ell \cdot 1}$ denote $q_{1,\ldots, 1}$ and similarly let
$q_{k\cdot 0, \ell\cdot 1}$ denote $q_{0,  \ldots, 0, 1, \ldots, 1} = {q_{1,\ldots, 1}}^{*2^{k}}$
in $H_{2^{k}(2^{\ell} - 1)}(B \si_{2^{k + \ell}})$, where there are $k$ zeros and $\ell$ ones.  

Let $\gamma_{\ell, n}$ denote the linear dual to ${(q_{\ell \cdot 1})}^{*n}$
in the Nakaoka monomial basis.
\end{definition}

If $\alpha \in H^{*}(\coprod_{n} B \si_{n})$ is a monomial in the $q_{I}$ we 
let $\alpha^{\vee} \in 
H^{*}(\coprod_{n} B \si_{n})$ denote the cohomology class which evaluates to one
on $\alpha$ and is zero on all other monomials. 
In particular $\gamma_{\ell, 2^{k}}$ is  
${q_{k \cdot 0, \ell \cdot 1}}^{\vee}$.  We will use these $\gamma_{\ell, n}$ as Hopf ring generators for the cohomology of symmetric groups, thus making the following the first  step in geometrically representing this cohomology.

\begin{definition}
Let $\Gamma_{\ell,n}$ be
defined as the collection of ${\bf x} = (x_{1}, \cdots, x_{n2^{\ell}})/\sim$
which can be partitioned into $n$ sets of $2^{\ell}$ points such that all points in each set share their first coordinate. 
\end{definition}

$\Gamma_{\ell,n}$ is the image in $\uconf_{n\cdot 2^{\ell}}(\R^{d})$ 
of $\widehat{\Gamma_{\ell,n}},$ in which along with the points there is a choice of partition.

\begin{theorem}\label{T:geomrep}
The cohomology class $\gamma_{\ell, n}$ is the Thom class of the variety $\Gamma_{\ell,n}$.
\end{theorem}

We first record the following, which is immediate algebraically from the definition.

\begin{lemma}\label{L:AlgCoprod}
The coproduct of $\gamma_{\ell, n}$ is given by
$$\Delta \gamma_{\ell, n} = \sum_{i+j = n} {\gamma_{\ell, i}} \otimes {\gamma_{\ell, j}}.$$
\end{lemma}

\begin{figure}
$$\includegraphics[width=4cm]{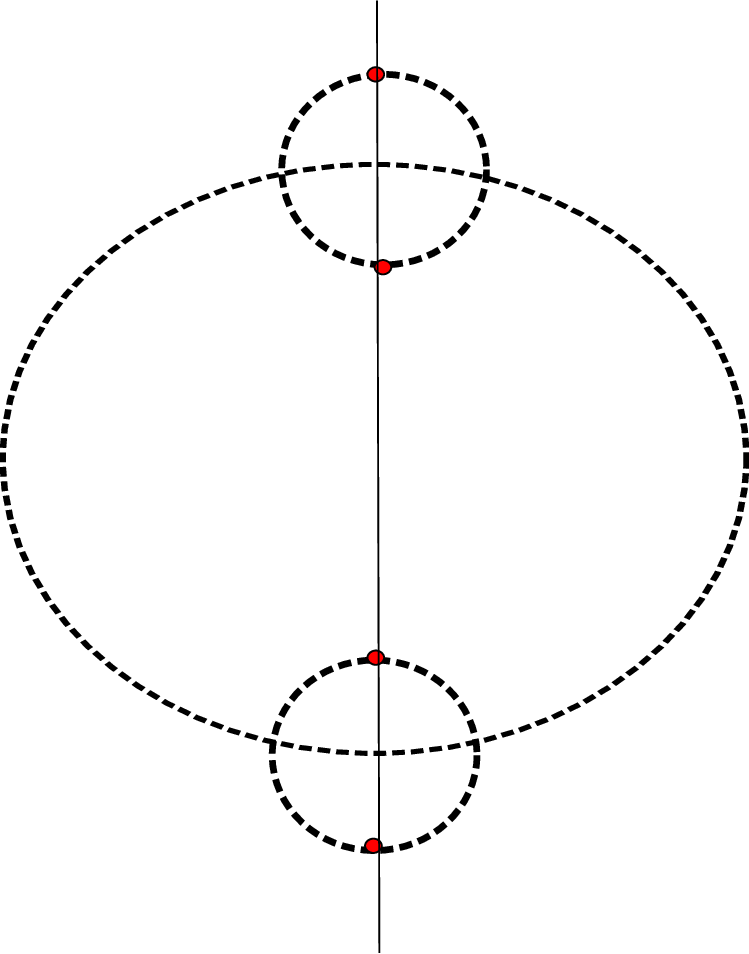}$$
\caption{An illustration that $\Gamma_{2,1}$ intersects $\orb_{1,1}$ exactly once.}\label{fig2}
\end{figure}

\begin{proof}[Proof of \refT{geomrep}]
We start by showing that 
the coproduct formula of \refL{AlgCoprod} holds for the the Thom class of $\Gamma_{\ell, n}$.  We model 
the product map by the embedding $\uconf_{i 2^{\ell}}(\R^{\infty})  \times \uconf_{(n-i)2^{\ell}}(\R^{\infty})  \to 
\uconf_{n2^{\ell}}(\R^{\infty})$ by using homeomorphism of $\R$ with the negative (respectively positive) 
real numbers to change the first coordinates of the first (respectively second) given configurations and taking
their union.  This  model of the product map is transversal to $\Gamma_{\ell, n}$, whose preimage in 
$\uconf_{i 2^{\ell}}(\R^{\infty})  \times \uconf_{(n-i)2^{\ell}}(\R^{\infty})$ is exactly 
$\Gamma_{\ell, i} \times \Gamma_{\ell, n-i}.$  Because the Thom class of a preimage of a subvariety under a transversal
map is the pull-back of its Thom class, the desired coproduct formula follows.

By Nakaoka's calculation as stated in \refT{Nakaoka}, the indecomposables in homology of $\bigsqcup_{n} B\si_{n}$
under the product lie in dimensions greater than $n - 1$ on components indexed by $n$ which are powers of two.  The product
map is thus surjective in lower degrees, or dually the coproduct map is injective, which implies that we may
use these coproduct formulae inductively to reduce to showing that the Thom class of $\Gamma_{\ell, 1}$ is
$\gamma_{\ell, 1}$.   Still using the coproduct formula, along with compatibility of evaluation of cohomology on  homology
with the product and coproduct, the value of $\Gamma_{\ell, 1}$ on any product is zero.
In degree $(2^{\ell} - 1)$ the only indecomposible in homology is $q_{\ell \cdot 1}$.  We can 
check immediately that $\orb_{\ell \cdot 1}$ which represents $q_{\ell \cdot 1}$
intersects with $\Gamma_{\ell, 1}$ in exactly one point, and the tangent vectors span the full tangent space of the configuration space as needed for transversality, as we illustrate in Figure~\ref{fig2}.
\end{proof}

These varieties $\Gamma_{\ell, n}$ are analogues of Schubert varieties, as we will see more precisely
in Section~\ref{S:SW}.
We can use other coordinates, or codimension one subspaces, to define $\Gamma_{\ell,n}$ and then use
the geometry of cup and transfer products to understand  representing varieties.  For example,
$\gamma_{1,2} \gamma_{2,1}$ is Thom class of the subvariety defined 
by ``four points which share their first coordinate and break up into
two groups of two points which share their second coordinate,'' while $\gamma_{1,1}^{3} \tr \gamma_{1,1}$ is
the Thom class of the subvariety defined  by configurations with
``two points which share their first three coordinates and another two which share their fourth 
coordinate.''

\bigskip

Getting back to homology, 
on the $q_{I}$ the coproduct dual to the cup product is classically known, and thus it is determined on the entire homology of symmetric groups because of the bialgebra structure.

\begin{definition}\label{D:diagonal}
Define a coproduct $\Delta_{\cupp}$ on $H_{*}(\coprod_{n} B \si_{n})$ by extending the formula for $I$ admissible $\Delta_{\cupp} (q_{I}) = \sum_{J+K = I} q_{J} \otimes q_{K}$ , where when $I = i_{1}, \cdots, i_{n}$ we have that $J$ and $K$ range over partitions of the same length such that for each $\ell$, $j_{\ell} + k_{\ell} = i_{\ell}$.
\end{definition}

This coproduct is more complicated than it seems at first.  Even when starting with an admissible $I$, the sum above is over all possible $J$ and $K$.  Thus to get an expression in the standard basis, as
needed for example to apply the coproduct again, one must apply Adem and $q_{0}$ relations.  The ones which get used most often are the relations $q_{2n+ 1}  q_{0} = 0$ and $q_{2n}  q_{0} =  q_{0} q_{n}$. 

\begin{theorem}[See for example I.2 of \cite{CLM76}]\label{T:cupcoprod}
Under the isomorphism of \refT{Nakaoka}, the diagonal map on 
$\coprod_{n} B \si_{n}$ induces the map $\Delta_{\cupp}$ on homology.
\end{theorem}

On the other hand, one of our main results is that the coproduct dual to the transfer product is primitive.
We originally proved this geometrically, but now by Kuhn's suggestion we use the machinery of \cite{BMMS86}.

\begin{theorem}\label{T:primitive}
The transfer product is linearly dual to the primitive coproduct on the Kudo-Araki-Dyer-Lashof algebra.  That is, $\Delta_{\tr}(q_{I}) = q_{I} \otimes 1 + 1 \otimes q_{I}$, where $1$ is the non-zero class in $H_{0}(B \si_{0})$.
\end{theorem}

\begin{proof}
Let $D_k(Y)=E\si_k \ltimes_{\si_k} Y^k$, for a based space $Y$.
We recall that $D(Y)=\bigvee_{k \geq 0} D_k(Y)$ 
is the free $E_\infty$-space $Y$, so in particular  $D(S^0) = \bigvee_{n \geq 0}B\si_n$. 
Similarly for a spectrum $\si$ we denote 
by $\mD(\si)$ the free $E_\infty$-spectrum it generates.
In particular $\mD(\Sigma^\infty Y) \simeq \Sigma^\infty D(Y)$.


As mentioned in Remark~\ref{stable}, the transfer product is
 induced by a stable map  
 $$\Psi: \Sigma^\infty D(S^0) \to \Sigma^\infty \left( D(S^0) \sm D(S^0) \right).$$  
 Observe that there is an equivalence
$\beta:D(S^0 \vee S^0) \simeq D(S^0) \sm D(S^0)$.

Recall Theorem~\ref{T:Nakaoka}, which says that
$H_*(D(S^0))$ is the free $\K$-algebra on the generator $\iota
\in \tilde{H}_0(S^0)$.
Similarly,   $H_*(D(S^0 \vee S^0))$ is the free $\K$-algebra 
on the two generators 
$$	(\iota,0),(0,\iota) \in \tilde{H}_0(S^0 \vee S^0).$$
 Clearly $\beta_*(\iota,0)=\iota \otimes 1$ and 
$\beta_*(0,\iota)=1 \otimes \iota$. 

Theorem 1.5 in \cite{BMMS86} states that 
$\Psi$ can be identified to $\mathcal{D}(v):\mD(\und{S}^0) \to \mD(\und{S}^0 \vee \und{S}^0)$, where $v:\und{S}^0 \to \und{S}^0 \vee \und{S}^0$ is the pinch map of the sphere spectrum. 

This implies that $\Delta_\tr$ is a map of $\K$-algebras. Explicitly, by the external Cartan formula,
$$\Delta_\tr q_n(a)= \sum_{i+j=n}  (q_i \otimes q_j) \Delta_\tr(a),$$
for $a \in H_*(\coprod_{n} B \si_{n})$.
The class $\iota$ is clearly primitive since $v_*(\iota)=(\iota,0)+(0,\iota)$ implies 
$\Delta_\tr(\iota)= \iota \otimes 1 + 1 \otimes \iota$.  
The external Cartan formula and the vanishing property of Definition~\ref{defq} 
imply that $q_I$ is primitive for each $I$ by induction on the length of $I$, completing the proof.

\end{proof}

We use this theorem to quickly determine
the cohomology of symmetric groups as a Hopf ring.

\section{Hopf ring structure through generators and relations}

The primitivity of the transfer coproduct  coupled with some 
classical theorems immediately leads to algebraic
presentations of $H^{*}(\coprod_{n} B \si_{n})$.
Recall from Theorem~7.15 of Milnor and Moore's standard reference \cite{MiMo65} that a Hopf algebra  which is polynomial and primitively generated has a linear dual that is exterior, generated by linear duals (in the monomial basis) to generators raised to powers of two.  
\refT{primitive} implies the following.

\begin{corollary}\label{C:exterior}
Under the transfer product $\tr$ alone, the cohomology of $\coprod_{n} B \si_{n}$ is an exterior algebra, generated by $({q_{I}}^{*2^{k}})^{\vee}$ for $I$ strongly admissible, or equivalently by ${q_{I}}^{\vee}$ for $I$  admissible. 
\end{corollary}

To incorporate the cup product structure, we appeal to another classical theorem.
As in I.3 of \cite{CLM76}, let $R[n]$ be the span of the $q_{I}$ of length $n$, a submodule of $H_{*}(B \si_{2^{n}})$.  
 Recall that $\gamma_{\ell, m}$ is the linear dual to 
${(q_{\ell \cdot 1})}^{*m}$ in the Nakaoka basis.

\begin{theorem}[Theorem~I.3.7 of \cite{CLM76}]\label{T:Rn}
The linear dual of $R[n]$, which is an algebra under cup product, is isomorphic to the polynomial algebra
generated by the classes $\gamma_{\ell, 2^{k}}$ with $k + \ell = n$, which we denote $D_{n}$.
\end{theorem}

In Section~\ref{S:invt} we review the classical fact that $D_{n}$ is canonically isomorphic to the $n$th Dickson algebra.
For a sketch of proof of \refT{Rn}, it is simple to see that  $q_{0, \ldots, 0, 1, \ldots, 1}$ are primitive under $\Delta_{\cupp}$
because $q_{1} \circ q_{0} = 0$.  It is a straightforward  induction to show there are no other primitives,
and then a counting argument to show that there are no relations among the 
${q_{k\cdot0, \ell \cdot 1}}^{\vee}$.  Note however that because of the Adem relations
in the Kudo-Araki-Dyer-Lashof algebra, the pairing between $q_{I}$ and polynomials
in $\gamma_{\ell, 2^{k}} = (q_{\ell \cdot 0, k \cdot 1})^{\vee}$ is complicated.   
For example, ${\gamma_{1,2}}^{3} = ({q_{0,1}}^{\vee})^{3} = {q_{0,3}}^{\vee} + 
{q_{2,2}}^{\vee},$ in part since $\Delta_{\cupp} q_{2,2}$ includes as a term
$q_{2,0} \otimes q_{0,2}$, which is equal to $q_{0,1} \otimes q_{0,2}$.

\refT{Rn} gives us the last input we need to understand the
cohomology of symmetric groups as a Hopf ring, since we see that under $\tr$ alone
the generators are polynomials in the $\gamma_{\ell, 2^{k}}$.

\begin{theorem}\label{T:AlgMain}
As a Hopf ring, $H^{*}(\coprod_{n} B \si_{n})$ is generated by the classes 
$\gamma_{\ell, 2^{k}}$.

The transfer product is exterior, and the antipode map is the identity.  
The $\gamma_{\ell, 2^{k}}$ with $\ell + k = n$ form a polynomial
ring.
\end{theorem}

The stated facts along with Hopf ring distributivity and the fact that the products of classes on different
components are zero determine the cup product structure, in particular of individual components.

\begin{corollary}\label{C:GenCond}
Any collection of classes $\{ \alpha_{\ell, 2^{k}} \}$ such that 
$ \alpha_{\ell, 2^{k}} \in H^{2^k(2^{\ell}-1)} (B \si_{2^{k + \ell}})$ pairs non-trivially with 
$q_{k\cdot 0, \ell \cdot 1}$ constitutes a generating set for $H^{*}(\coprod_{n} B \si_{n})$ 
as a Hopf ring.
\end{corollary}

In order to apply
the distributivity relation, we need the coproduct
of $\gamma_{\ell, 2^{k}}$ as given by \refL{AlgCoprod}.
Finally, we compute transfer products by taking the binary expansion 
$ j = \sum 2^{k_{i}}$ with $k_{i}$ distinct, so that ${q_{\ell \cdot 1}}^{*j} =
 \prod {q_{\ell \cdot 1}}^{*2^{k_{i}}}$.  By \refT{primitive} and linear duality it follows that
\begin{equation*}
\gamma_{\ell, j} 
= \bigodot_{j = \sum 2^{k_{i}}} \gamma_{\ell, 2^{k_{i}}}.
\end{equation*}
Because the transfer product is exterior, we obtain the following.

\begin{proposition}\label{P:p1}
The transfer products of classes $\gamma_{\ell,n}$ are
given by 
$$\gamma_{\ell, n} \tr \gamma_{\ell, m} = \binom{n+m}{n} \gamma_{\ell, n+m},$$
while transfer products between other classes have no relations.
\end{proposition}

Collecting \refT{AlgMain}, \refL{AlgCoprod}, \refT{Rn} and Proposition~\ref{P:p1} yields a
proof of \refT{genandrel}, our first presentation of the cohomology of symmetric
groups as a Hopf ring.

\refT{AlgMain} and the fact that transfer product is exterior 
leads to an additive basis for the cohomology of symmetric groups.
Since the $\gamma_{\ell, 2^{k}}$ are Hopf ring generators, their Hopf ring monomials span, and an induction
by number of transfer products using \refT{primitive} shows they are independent.  Let ${\mathcal{P}}(m)$ denote the set of partitions of $m$ into nonnegative powers of two.  If $P$ is such a partition, let $P_{n}$ denote the number of times $2^{n}$ occurs, and let $D_{P}$ denote 
$\bigotimes_{n} \bigwedge^{P_{n}} D_{n}$ where $D_{n}$ is
the polynomial algebra on $\gamma_{\ell, 2^{k}}$ with $\ell + k = n$ as stated in \refT{Rn}.  By convention $D_{0}$
is the ground field.  We recover an additive
isomorphism well-known to experts.

\begin{proposition}
As a graded vector space, $H^{*}(BS_{n})$ is isomorphic to $\bigoplus_{P \in {\mathcal{P}}(n)} D_{P}$.
\end{proposition}

It is straightforward to calculate the Poincar\'e polynomial for this cohomology.  We have found it
unenlightening since we could only  describe Poincar\'e polynomials of exterior powers
of polynomial algebras using inclusion-exclusion methods. 

In the next section we develop a slightly improved basis, which requires the following.

\begin{proposition}\label{P:p2}
The classes $\{ \gamma_{\ell, n} \}$ such that $n 2^{\ell} = m$ generate a polynomial
subalgebra of $H^{*}(B \si_{m})$.
\end{proposition}

\begin{proof}[Sketch of proof.]
We start with $m = 2^{p}$ which is covered by 
\refT{Rn}, which says that the classes 
$\gamma_{\ell, 2^{k}}$
with $k + \ell = p$ form a polynomial subring of $H^{*}(B\si_{2^{p}})$.  
We then use in induction on the number of ones in the binomial expansion
of $m$ and Hopf ring distributivity to establish the general case.
\end{proof}

\begin{figure}
\begin{tabular}{c|cccccccc}
15&&&&&&&&$\gamma_{4,1}$\\
14&&&&&&&&$\gamma_{3,2}$\\
13&&&&\\
12&&&&&&&&$\gamma_{2,4}$\\
11&&&&\\
10&&\\
9&&&&&&$\gamma_{2,3}$\\
8&&&&&&&&$\gamma_{1,8}$\\
7&&&&$\gamma_{3,1}$&&&$\gamma_{1,7}$\\
6&&&&$\gamma_{2,2}$&&$\gamma_{1,6}$\\
5&&&&&$\gamma_{1,5}$\\
4&&&&$\gamma_{1,4}$\\
3&&$\gamma_{2,1}$&$\gamma_{1,3}$\\
2&&$\gamma_{1,2}$\\
1&$\gamma_{1,1}$\\
\hline
&$B\si_2$&$B\si_4$&$B\si_6$&$B\si_8$&$B\si_{10}$&$B\si_{12}$&$B\si_{14}$&$B\si_{16}$
\end{tabular}

\caption{Hopf ring generators of  $H^*(\coprod_{n} B \si_{n})$ through $B\si_{16}$. }
\end{figure}


By \refT{AlgMain} that $\gamma_{\ell, n}$ are Hopf ring generators, \refT{geomrep} that these classes
are Thom classes of linear subvarieties, and the 
geometric interpretations of cup and transfer products,  all of the cohomology of symmetric groups
is represented by such subvarieties.  These are defined by groupings and subgroupings 
of  points into sets with cardinalities which are powers of two which share coordinates.
The third author began this investigation by conjecturing that representations by linear subvarieties
would be possible, since
linear submanifolds represent classes in ordered configuration spaces.

\section{Presentation of product structures through an additive basis}\label{S:comps}

We can use our knowledge of the Hopf ring structure on $H^{*}(\coprod_{n} B \si_{n})$ 
to explicitly understand the cup and transfer product structures
through additive bases.  We set the notational conventions that cup product 
has priority over transfer product, so that $a \cupp b \tr c$ means $(a \cupp b) \tr c$, and that
exponents always refer to repeated application of cup product (an easy choice, since transfer
product is exterior).  Let $1_{m}$ denote
the unit for cup product on component $m$.  

We proceed with some calculations on the first even components 
(since the $\F_{2}$-cohomology
of $B \si_{2k + 1}$ is isomorphic to that of $B \si_{2k}$.)
As $B \si_{2} \simeq \R P^{\infty}$  its cohomology is a polynomial ring generated
by $\gamma_{1,1}$.  

For $H^{*}(B \si_{4})$ the Hopf ring monomial basis
 consists of classes ${\gamma_{1,1}}^{i} \tr {\gamma_{1,1}}^{j} \in H^{i+j}(B \si_{4})$ (which
 are zero if $i=j$) along
with polynomials in $\gamma_{1,2} \in H^{2}(B \si_{4})$ and 
$\gamma_{2,1} \in H^{3}(B \si_{4})$.  
Using Hopf ring distributivity, we have that
$$({\gamma_{1,1}}^{i} \tr {\gamma_{1,1}}^{j}) \cupp ({\gamma_{1,1}}^{k} \tr {\gamma_{1,1}}^{\ell})
= {\gamma_{1,1}}^{i+k} \tr {\gamma_{1,1}}^{j+ \ell} + {\gamma_{1,1}}^{i+\ell} \tr {\gamma_{1,1}}^{j+k},
$$
where we note one of these terms could be zero, if either $i+k = j+ \ell$ or
if $i + \ell =  j + k$.
In order to compute some products with $\gamma_{1,2}$ we have to use its coproduct,
which by \refL{AlgCoprod} is equal to 
$\gamma_{1,2} \otimes 1_{0} + \gamma_{1,1} \otimes \gamma_{1,1} + 1_{0} \otimes  \gamma_{1,2}$.
Using distributivity,
\begin{multline*}
\gamma_{1,2}\cupp ({\gamma_{1,1}}^{n} \tr {\gamma_{1,1}}^{m}) = 
(\gamma_{1,2} \cupp {\gamma_{1,1}}^{n}) \tr (1_{0} \cupp
{\gamma_{1,1}}^{m}) +\\  (\gamma_{1,1} \cupp {\gamma_{1,1}}^{n}) \tr (\gamma_{1,1} \cupp
{\gamma_{1,1}}^{m}) +  (1_{0} \cupp {\gamma_{1,1}}^{n}) \tr (\gamma_{1,2} \cupp
{\gamma_{1,1}}^{m}) = {\gamma_{1,1}}^{n+1} \tr {\gamma_{1,1}}^{m+1}.$$
\end{multline*}
In general, most terms arising from the Hopf ring distributivity relation are zero
because they involve multiplication of classes supported on different components.
The last basic products to compute for $B \si_{4}$
are  $\gamma_{2,1}\cupp ({\gamma_{1,1}}^{n} \tr {\gamma_{1,1}}^{m})$,
which are zero because the coproduct of $\gamma_{2,1}$ is just 
$\gamma_{2,1} \otimes 1_{0} + 1_{0} \otimes \gamma_{2,1}.$ 
Applying distributivity repeatedly  we get that if $k \neq 0$
$${\gamma_{1,2}}^{p} {\gamma_{2,1}}^{q} \cupp ({\gamma_{1,1}}^{k} \tr {\gamma_{1,1}}^{\ell})
= 
\begin{cases}
{\gamma_{1,1}}^{k+p} \tr {\gamma_{1,1}}^{\ell+p}  \;\;\; & {\rm if} \;\; q = 0.\\
0 \;\;\;\;\;\;  & {\rm if} \;\; q \neq 0,
\end{cases}
$$
which completes an understanding of how to multiply elements of our additive basis.
The case of $B \si_{4}$ is one of the very few in which it is simpler to understand 
the cup multiplicative structure in terms of ring generators and relations.
From the multiplicative rules just given, it is a straightforward exercise to 
deduce that $\gamma_{1,1} \tr 1_{2}$, 
$\gamma_{1,2}$ and $\gamma_{2,1}$ generate the cohomology
on this component, with the lone relation being $(\gamma_{1,1} \tr 1_{2}) \cupp \gamma_{2,1} = 0$.
Similarly, one can write down an additive basis for 
$H^{*}(B \si_{6})$, determine its multiplication rules, and then show that it is 
generated by $\gamma_{1,1} \tr 1_{4}$, $\gamma_{1,2} \tr 1_{2}$, $\gamma_{2,1} \tr 1_{2}$ and 
$\gamma_{1,1}^{2} \tr \gamma_{1,1} \tr 1_{2}$, with the relation that $\gamma_{2,1} \cupp 
({\gamma_{1,1}}^{2} \tr \gamma_{1,1} \tr 1_{2}) = 0$, in agreement with the results of Chapter VI.5 of 
\cite{AdMi94}.  

\begin{remark}
We can also see relations through our geometric representatives for cohomology.  For example,
$ \gamma_{2,1} \cupp (\gamma_{1,1} \tr 1_{2}) $ is represented by the subvariety of ``four points which
share their first coordinate, two of which share their second coordinate.''  This subvariety is cobounded by
``four points, two of which share their first coordinate, two of which share their first and second coordinate,
with the first two having a first coordinate which is less than that of the second two.''  
\end{remark}

\medskip

In general, presentations in terms of generators and relations are quite complicated.
We instead understand cup and transfer products 
explicitly in terms of a canonical additive basis.
Recall the notion of Hopf ring monomial from \refD{Hopfmono}.

\begin{definition}
A gathered monomial in the cohomology of symmetric groups is a
Hopf ring monomial  in the generators $\gamma_{\ell, n}$  
where such $n$ are maximal or equivalently the number of transfer products
which appear is minimal.
\end{definition}

For example, $\gamma_{1,4} {\gamma_{2,2}}^{3} \tr \gamma_{1,2} {\gamma_{2,1}}^{3}
= \gamma_{1,6}  {\gamma_{2,3}}^{3}$.  Gathered monomials such as the latter in which 
no transfer products appear are building blocks for general gathered monomials.

\begin{definition}
A gathered block is a monomial of the form $\prod_{i} {\gamma_{\ell_{i}, n_{i}}}^{d_{i}},$
where the product is the cup product.  Its profile is defined to be the collection
of pairs $(\ell_{i}, d_{i})$.

Non-trivial gathered blocks must have all of the numbers $2^{\ell_{i}} n_{i}$
equal, and we call this number divided by two the width.  
We assume that the factors are ordered from smallest to largest $n_{i}$ (or largest
to smallest $\ell_{i}$), and then note that $n_{i} = 2^{\ell_{1} - \ell_{i}} n_{1}$.
\end{definition}

\begin{proposition}
A gathered monomial can be written uniquely as the transfer product of gathered blocks with distinct profiles. 
Gathered monomials form a canonical additive basis for the cohomology of $\coprod_{n} B \si_{n}$.
\end{proposition}

\begin{figure}
\centering
\subfloat[$(\gamma_{1,1}^i \tr \gamma_{1,1}^j) \cupp (\gamma_{1,1}\tr 1_2) =\gamma_{1,1}^{i+1} \tr \gamma_{1,1}^j+\gamma_{1,1}^{i} \tr \gamma_{1,1}^{j+1}$]{
\hspace{1cm}
\begin{tabular}{ccccccc}
 \includegraphics[height=48pt]{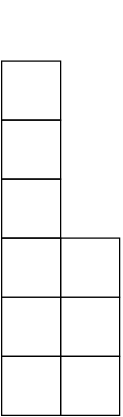}& \raisebox{18pt}{$\cupp$} &\includegraphics[height=48pt]{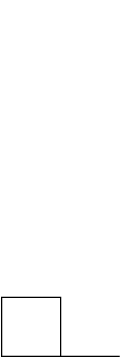} &\raisebox{18pt}{=}&  \includegraphics[height=48pt]{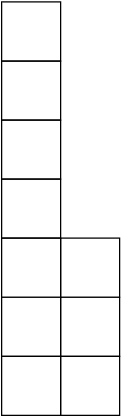} &\raisebox{18pt}{+}&  \includegraphics[height=48pt]{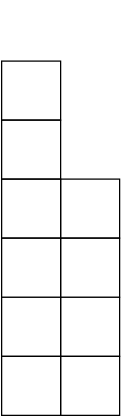}\\
\end{tabular}\hspace{1cm}}

\subfloat[$(\gamma_{1,1}^i \tr \gamma_{1,1}^j) \cupp \gamma_{1,2}=\gamma_{1,1}^{i+1} \tr \gamma_{1,1}^{j+1}$]{
\hspace{.5cm}
\begin{tabular}{ccccc}
 \includegraphics[height=48pt]{Figures/diag_mult_ij.eps}& \raisebox{18pt}{$\cupp$} &\includegraphics[height=48pt]{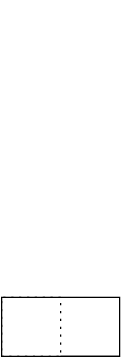} &\raisebox{18pt}{=}&  \includegraphics[height=48pt]{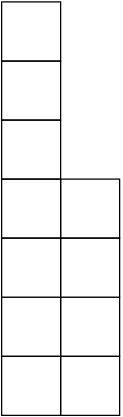}\\
\end{tabular}\hspace{.5cm}}
\qquad
\subfloat[$(\gamma_{1,1}^i \tr \gamma_{1,1}^j) \cupp\gamma_{2,1}=0$]{
\begin{tabular}{ccccc}
 \includegraphics[height=48pt]{Figures/diag_mult_ij.eps}& \raisebox{18pt}{$\cupp$} &\includegraphics[height=48pt]{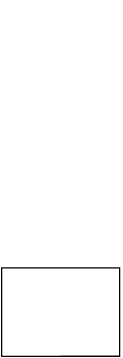} &\raisebox{18pt}{=}&  \raisebox{18pt}{0}\\
\end{tabular}
}

\subfloat[$(\gamma_{1,1}^3 \tr \gamma_{2,1}\tr 1_2) \cupp (\gamma_{1,2} \tr 1_{4}) =\gamma_{1,1}^4\tr\gamma_{2,1}\tr 1_2 + \gamma_{1,1}^3\tr\gamma_{1,2}\gamma_{2,1}\tr 1_2$]{
\hspace{1cm}

\begin{tabular}{ccccccc}
 \includegraphics[height=32pt]{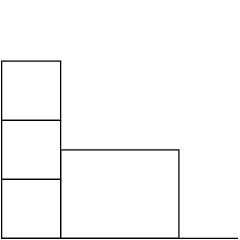}& \raisebox{14pt}{$\cupp$} &\includegraphics[height=32pt]{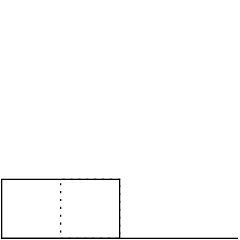} &\raisebox{14pt}{=}&  \includegraphics[height=32pt]{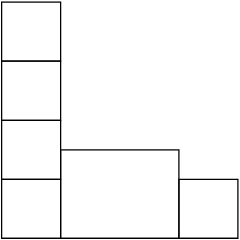}&\raisebox{14pt}{+}&\includegraphics[height=32pt]{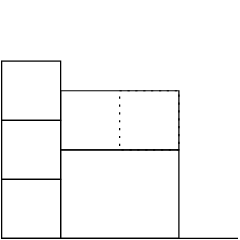}\\
\end{tabular}\hspace{1cm}

}
\caption{Some computations in $H^*(B\si_{4})$ and $H^*(B\si_{8})$, expressed by both gathered monomials and skyline diagrams.}
\label{F:skyline}
\end{figure}

Representing gathered monomials graphically is helpful.
We represent $\gamma_{\ell,n}$ by a rectangle of width $n \cdot 2^{\ell}$ and height $1 - \frac{1}{2^{\ell}}$, so that its 
area corresponds to its degree. 
We  represent $1_{n}$ by an edge of width $n$ (a height-zero rectangle).
A gathered block, which is a product of $\gamma_{\ell, n}$ for fixed $n \cdot 2^{\ell}$,
is represented by a single column of such rectangles, stacked on top of each other, with order which does not matter.
A gathered monomial is represented by placing such columns next to each other, which we call the skyline diagram of the monomial.  We also refer to the gathered monomial basis as the skyline basis to emphasize
this presentation.
See Figure~\ref{F:skyline} below for an illustration.

\begin{definition}
Let  $\prod_{i} {\gamma_{\ell_{i}, n_{i}}}^{d_{i}}$ be a gathered block, and 
let $n_{1} = \sum_{j = 1}^{k} m_{j}$ be a partition of $n_{1}$.
A partition of this gathered block into $k$  is defined by the set consisting of the $k$ blocks
$ \prod_{i} {\gamma_{\ell_{i}, 2^{\ell_{1} - \ell_{i}}m_{j}}}^{d_{i}}.$  We allow for some $m_{j}$
to be zero, in which case the corresponding elements of the partition will be $1_{0}$.

A splitting of a gathered monomial $f_{1} \tr \cdots \tr f_{k}$ into two
 is a pair of gathered monomials $f_{1}' \tr \cdots \tr f_{k}' $
and $f_{1}'' \tr \cdots \tr f_{k}''$ where each $\{ f_{i}', f_{i}'' \}$ is a
partition of $f_{i}$ into two  (which could be trivial - that is, of the form
$\{ 1_{0}, f_{i} \}$).
\end{definition}

\begin{proposition}\label{P:gatheredcoprod}
The coproduct of a gathered monomial is given by 
$$\Delta  f_{1} \tr \cdots \tr f_{k} = \left( \sum f_{1}' \tr \cdots \tr f_{k}' \right) \otimes \left( f_{1}'' \tr \cdots \tr f_{k}''\right),$$ 
where the sum is over all splittings of the monomial into two.
 \end{proposition}

\begin{proof}
To establish the special case of gathered blocks - that is, having only one $f$ - 
we use the Hopf algebra compatibility of cup product and Pontryagin coproduct.  The 
coproduct $\Delta$ of any  ${\gamma_{\ell_{i}, n_{i}}}^{d_{i}}$ will correspond to partitions
of $n_{i}$ into two.  But only for the partitions of $n_{1}$ will there
be  corresponding 
partitions for all $n_{i}$ which yield non-trivial classes when cupped together.  The 
resulting products correspond to the partitions of $f$ into two.

The general case follows from the Hopf algebra compatibility of partition product and
Pontryagin coproduct.  Because the monomial is gathered, no terms in the coproducts
of $f_{i}$ can be equal, so we obtain no trivial transfer products when such terms are collected.
\end{proof}

In terms of skyline diagrams, the coproduct can be understood by introducing vertical dashed lines in the
rectangles representing $\gamma_{\ell,n}$, dividing the rectangle into $n$ equal pieces.  The coproduct
is then given by dividing along all existing columns and vertical dashed lines of full height
and then partitioning them into two to make two new skyline diagrams.

\begin{definition}
A partition of a gathered monomial in $H^{*}(B \si_{m})$
is a partition of each of its gathered blocks.
The associated component partition is the partition of $m$ given by the components of the classes
in the partition.

We define the refinement of a partition of a gathered monomial in the obvious way,
reflected faithfully by the
refinement structure of the associated component partitions.

A matching $\mu$ between partitions of two gathered monomials is an isomorphism of 
their respective component partitions.  We say that one matching refines another 
if that isomorphism commutes with inclusions of components under refinement. 
\end{definition}

For any gathered monomial in $H^{*}(B\si_{m})$  there is a canonical partition of
$m$ defined by the components of its constituent gathered block monomials.  
The associated component
partition of a monomial partition is a  refinement of this canonical partition.

We are now ready to describe product structures in terms of our additive basis
of gathered monomials.

\begin{theorem}\label{T:generalcup}
The transfer product of two gathered monomials ${\bf m}$ and ${\bf n}$
is a multiple of the gathered monomial
whose gathered block of a given profile has width which is the sum of the widths of 
the blocks of that profile in ${\bf m}$ and ${\bf n}$.
The multiple is zero if and only if any of those two widths share some non-zero digit
of their binary expansion.  


Let $M_{{\bf m},{\bf n}}$ denote 
the set of matchings between any of the partitions of these gathered monomials
which are not  a refinement of some other matching.
The cup product of $x$ and $y$ is the sum 
$$\sum_{\mu \in M_{{\bf m},{\bf n}}} \left( \bigodot_{b, b' {\textrm{ matched by }} \mu}  \beta_{\mu} b \cupp b'
\right),$$
where $\beta_{\mu}$ is zero if and only if there are two products $b \cupp b'$ which result in blocks with
the same profile and whose widths have binary expansion which share a non-zero digit.
\end{theorem}

Graphically, transfer product corresponds to placing two column Skyline diagrams next to each
other and merging columns with the same constituent blocks, with a coefficient of zero if any of those column
widths share a one in their dyadic expansion.  
For cup product, we start with two column diagrams and consider all possible ways to 
split each into columns, along either original boundaries of columns or along the vertical lines of full height
internal to the rectangles representing $\gamma_{\ell,n}$.  We then match columns of each in 
all possible ways up to automorphism, and stack the resulting matched columns to get a new
set of columns -- see the  Figure~\ref{F:skyline}.

\begin{proof}[Proof of \refT{generalcup}]
We use gathered blocks, whose multiplication is polynomial by \refP{p2} 
as a base case for an induction on the total number of
blocks in ${\bf m}$ and ${\bf n}$.  View say ${\bf m}$
as a non-trivial transfer product of ${\bf m}'$ and ${\bf m}''$ which preserves blocks, 
so that ${\bf m}'$ and ${\bf m}''$ each has fewer blocks than ${\bf m}$.  
The key is to see that each matching in $M_{{\bf m},{\bf n}}$
coincides with some (arbitrary) partition of $y$ into $\{{\bf n}', {\bf n}''\}$ 
along  with matchings  of partitions of those pieces with partitions of ${\bf m}'$ and ${\bf m}''$.
From this observation, the induction follows, with the coefficient $\beta_{\mu}$ 
accounting for when such a process yields a partition product of some monomial
in the $\gamma_{\ell, 2^{k}}$ with itself.
\end{proof}

Given that the basis of skyline diagrams is a fundamental cohomology basis,
it would be helpful to understand the
pairing of gathered monomials with Nakaoka's monomial basis for homology.  Polynomials in 
$\gamma_{\ell, 2^{k}}$ are the fundamental case, which as mentioned after
\refT{Rn} pair non-trivially with the basis of $q_{I}$ with $I$ of length $k + \ell$.

\section{Topology and the invariant theoretic presentation}\label{S:invt}

Compare the presentation for the cohomology of symmetric groups as a 
Hopf ring, as given in \refT{genandrel}, with the Hopf ring presentations of rings of symmetric 
functions, as given in Example~\ref{Ex:classical} and \refP{symmfun}.  Seeing classes
which behave similarly, we obtain some immediate 
identifications of split quotient rings of the cohomology of symmetric groups.

\begin{definition}
Define the level-$\ell$ quotient Hopf ring of the cohomology 
of symmetric groups, denoted $\pa_{\ell}$, to be the 
quotient Hopf ring obtained
by setting all $\gamma_{\ell', n}$ for $\ell' \neq \ell$ equal to zero.
This quotient map is split by  the sub-Hopf ring generated
by the classes $\gamma_{\ell, n}$. 
Let $\pa_{\ell}[m]$ be the sub-module of $\pa_{\ell}$ supported on 
$B \si_{m}$, which is an algebra under cup product.
\end{definition}

Graphically, these sub-rings each consist of all skyline diagrams made from the blocks of one fixed size.

\begin{proposition}
The level-$\ell$ Hopf ring $\pa_{\ell}$ is isomorphic as Hopf ring to the total symmetric invariants of
a $\ka[x]$.  Thus $\pa_{\ell}[m]$ is a polynomial ring for any $m$.
\end{proposition}

The proof is an immediate comparison of their two presentations.  We originally
proved the second part directly from the Hopf ring presentation of $\pa_{\ell}[m]$,
before realizing that we were mimicking the proof that symmetric functions
form a polynomial algebra.

This identification has the following significant generalization.

\begin{definition}
Define the scale of $\gamma_{\ell, n}$ to be the product $\ell \cdot |n|_{2}$, where $|n|_{2}$ is 
the $2$-adic valuation of $n$ (that is, the largest power of two which divides $n$).
Define the scale-$k$ quotient of the cohomology of symmetric groups,
denoted $\scr_{k}$, to be the quotient Hopf ring obtained by setting all $\gamma_{\ell,n}$ with
either scale less than $k$ or with $\ell > k$ to zero.   It is isomorphic to the sub-Hopf ring
generated by $\gamma_{\ell,n}$ with scale greater than or equal to $k$ and $\ell \leq k$.
\end{definition}

Graphically, the skyline diagrams which are non-zero in this quotient are those 
made up of blocks of width exactly $2^{k-1}$.  

In the level-$\ell$ case, the split sub-ring was compatible with our additive basis, in that
if a sum of gathered monomials was in the sub-ring then each monomial was in the sub-ring
itself.  That is not true in the scale-$\ell$ setting, where for example ${\gamma_{1,2}}^{3}$
includes a term of ${\gamma_{1,1}}^{3} \tr {\gamma_{1,1}}^{2} \tr \gamma_{1,1}$.

\begin{proposition}\label{P:scale}
The scale-$k$ Hopf ring $\scr_{k}$ is isomorphic as a Hopf ring to the total symmetric invariants of a 
polynomial ring in $k$ variables.
\end{proposition}

Once again, the proof is by a direct comparison,  made possible by the
Hopf ring approach.  The canonical isomorphism between them
sends $\gamma_{\ell,n}$ with scale greater than or equal to $k$ and $\ell \leq k$  to the symmetric
polynomial $\sigma(\ell)_{m}$ with $m = {\frac{n}{2^{k-\ell}}}$.

Our goals in the rest of this section are twofold.  First we develop the standard topology which
underlies these isomorphisms.  Then, we move from 
these identifications of local invariant-theoertic sub/quotient rings to 
the global invariant-theoretic description of \refT{invt}.

\medskip

The predominant approach to the cohomology of symmetric groups has been through 
restricting cohomology to that of elementary abelian subgroups.  For the following, we refer to Chapters~3~and~4 in \cite{AdMi94}.
Let $V_{n}$ denote the subgroup of 
$(\Z/2)^{n} \subset \Sym_{2^{n}}$ defined by having $(\Z/2)^{n}$ act on itself.  
If we view this action
as given by linear translations on the $\F_{2}$-vector space 
$\oplus_{n}{\F_{2}}$, then we can see that the normalizer of this subgroup
is isomorphic to all affine transformations of $(\F_{2})^{n}$.  
The Weyl group is thus $GL_{n}(\F_{2})$,
which acts as expected on the cohomology of $V_{n}$.  The invariants 
$\F_{2}[x_{1}, \ldots, x_{n}]^{GL_{n}(\F_{2})}$ are known as Dickson algebras, which are
polynomial on generators $d_{k,\ell}$ in dimensions $2^{k}(2^{\ell} - 1)$ where $k + \ell = n$.
As mentioned earlier, these Dickson algebras together form a Hopf ring, which we are currently
investigating.

Since we base our work on Nakaoka's homology calculation, 
our analysis of elementary abelian subgroups involves homology as well as cohomology.

\begin{lemma}\label{L:img}
The image of the homology of $BV_{n}$ in that of $B\si_{2^{n}}$ is exactly the span of the $q_{I}$
for $I$ admissible of length $n$.
\end{lemma}

\begin{proof}
The inclusion of $V_{n}$ into $\si_{2^{n}}$ factors through the $n$-fold iterated
wreath product of $\Z/2$ with itself, that is
$\Z/2 \int \left( \Z/2 \int \left( \cdots \left( \Z/2 \int \Z/2 \right) \cdots \right) \right)$.   
A well-known alternate definition of the Dyer-Lashof operations $q_{i}$ is through
the homology of the inclusion of wreath products $\Z/2 \int \si_{n} \to \si_{2n}$.
Inductively, the image of this iterated wreath product 
is given by length-$n$ Kudo-Araki-Dyer-Lashof classes, so the image of $V_{n}$ is contained in
the span of such operations.

To see that the image of $V_{n}$ yields all such classes we compare ranks using the
dual map in cohomology.    We claim that the image in cohomology
of the inclusion of $V_{n}$ in $\si_{2^{n}}$ is all of the Dickson invariants
$\F_{2}[x_{1}, \ldots, x_{n}]^{GL_{n}(\F_{2})} \cong \F_{2}[d_{k, \ell}]$ with 
$k + \ell = n$ and $\ell > 0,$ a fact known by Milgram \cite{Milg09} which we share now.  
The standard representation of $\si_{2^{n}}$
through permutation matrices gives rise to a vector bundle.   Because 
$V_{n}$ embeds in $\si_{2^{n}}$ through the linear action of $(\F_{2})^{n}$ on itself,
on passing to a permutation representation the standard representation
yields the sum of all one-dimensional
real representations of $(\Z/2)^{n}$.
Thus when the corresponding bundle is pulled
back to $BV_{n}$ it splits as the sum of all possible line bundles.
So the total Stiefel-Whitney
class of this standard bundle in the cohomology of $B\si_{2^{n}}$ maps
to $\prod_{y \in H^{1}(BV_{n})} (1 + y)$, where $y$ ranges over linear combinations
of the $x_{i}$.  But classical invariant theory identifies
$\sum d_{k,\ell}$ with the product of all
$1 + \lambda$ where $\lambda$ varies over all linear functions in the $x_{i}$.
So these Stiefel-Whitney classes map exactly to the Dickson generators (or to zero).

By Madsen's calculation, Theorem~I.3.7 of \cite{CLM76} as recounted in \refT{Rn},  the linear dual to the span of
the length-$n$ $q_{I}$ is a polynomial algebra in classes of dimension $2^{k}(2^{\ell} - 1)$
with $k + \ell = n$ and $\ell > 0$.  Since the image of the cohomology of $B\si_{2^{n}}$ in that
of $BV_{n}$ has the same rank as this polynomial algebra, and thus as the span of $q_{I}$
of length $n$,
the image in homology must be all of this span.
\end{proof}

Because the only classes in the Weyl-invariant cohomology of $BV_{n}$ in degrees
$2^{k}(2^{\ell} - 1)$ are the Dickson classes, and the map in homology sends a generator
in that degree to $q_{k\cdot0, \ell \cdot 1}$, we have the following.

\begin{corollary}\label{C:maptodickson}
The restriction of $\gamma_{\ell,2^{k}}$ with $k + \ell = n$
to the elementary abelian subgroup $V_{n}$ is the Dickson class $d_{k, \ell}$.
\end{corollary}

Following Quillen \cite{Quil71, QuVe72},  Gunawardena-Lannes-Zarati \cite{GLZ88}
showed  that $H^{*}(B \Sym_{n})$ injects in the direct sum of the cohomology of elementary abelian
subgroups. 
\refL{img} gives an alternate proof for that theorem
through the following refinement.

\begin{corollary}\label{C:elem}
The image of the elementary abelian subgroup  $\prod_{i} V_{k_{j}}$
in the homology of any symmetric group which contains it (that is, of order $\sum 2^{k_{j}}$
or greater) is the span of products $\prod q_{I_{j}}$ where $I_{j}$ is of length $k_{j}$.
Thus, the map from the homology of all elementary abelian subgroups to the homology
of symmetric groups is surjective.
\end{corollary}

We now give a topological interpretation of \refP{scale}.

\begin{theorem}
The map from $H^{*}(\coprod_{n} B \si_{n})$ to its
image in the cohomology of $\coprod_{m} {BV_{k}}^{\frac{m}{2^{k}}}$
coincides with the quotient map defining the scale-$k$ quotient ring
$\scr_{k}$.
\end{theorem}

\begin{proof}
By \refC{elem}, the image in homology of $\coprod_{m} {BV_{k}}^{\frac{m}{2^{k}}}$ is
the submodule  of products of $q_{I}$ of length $k$.   By \refT{cupcoprod} and \refT{primitive}, it
is closed under the coproducts dual to cup and transfer product.  Thus the image of this map
of classifying spaces in cohomology, linear dual to this image of
homology, is a quotient of the cohomology of symmetric groups as a Hopf ring.  

Recall that
$\gamma_{\ell,n} = ({q_{\ell \cdot 1}}^{*n})^{\vee}$, so that all 
$\gamma_{\ell,n}$ with either scale less than $k$ or $\ell > k$ will evaluate to zero on
the image of homology.   An elementary counting argument shows that this ideal,
the quotient by which defines $\scr_{k}$, is as large as possible so that the restriction
of the cohomology of symmetric groups to these elementary abelian subgroups is
exactly $\scr_{k}$.
\end{proof}

We now give a global invariant theoretic description
of the cohomology of symmetric groups.

\begin{definition}
Consider the ring of polynomials $\F_{2}[x_{A}]$, where 
$A \subseteq \und{m} = \{ 1, \ldots, m\}$.  We call ${A'}$ a translate of ${A}$ if
they are disjoint and of the same cardinality.  A collection of translates is to be mutually disjoint.
Call a monomial $\prod x_{A_{i}}$ proper if whenever some $A_{i}$ and $A_{j}$ intersect,
one is contained in the other, say $A_{i} \subset A_{j}$, and $A_{j}$ is the union
of translates of $A_{i}$
\end{definition}

\begin{theorem}\label{T:invt}
The cohomology of
symmetric groups $H^{*}(\coprod_{n} B \si_{n}; \F_{2})$ is isomorphic to 
to the quotient of $\bigoplus_{m} \F_{2}[x_{A} | A \subset \und{m}]^{\si_{m}}$,
with $x_{A}$ in degree $2^{\#A} - 1$,
by the additive submodule consisting of symmetrizations of monomials which are
not proper.
\end{theorem}

\begin{proof}[Sketch of proof.]
We begin with the abstract Hopf ring description of $H^{*}(\coprod_{n} B \si_{n}; \F_{2})$
given in \refT{genandrel}, and show that it is isomorphic to the quotient stated.

Given some $A = \{1, \cdots, 2^{k}\} \subset \und{m}$ define its $i$th translate $\tau_{i} A$
to be $\{ 1 + i 2^{k}, \cdots, (i+1)2^{k}\}$.  
We start to define a map between $H^{*}(\coprod_{n}  B \si_{n}; \F_{2})$ and
this ring of invariants by sending 
$\gamma_{\ell, n}$ to $\prod_{i= 1}^{n} x_{\tau_{i} A}$ where $A = \{1, \cdots, 2^{\ell}\}$.
More generally, the gathered block $\prod {\gamma_{\ell_{i}, n_{i}}}^{d_{i}}$
maps to the symmetrized product where the first $n_{i}$ translates of $x_{1, \cdots, 2^{\ell_{i}}}$ are raised
to the $d_{i}$th power.  Such products are proper monomials. The transfer products
of gathered blocks go to symmetrized products of such monomials, after reindexing so that
the subscripts corresponding to different gathered blocks are distinct.   
Just as we showed for symmetric functions 
in \refP{symmfun}, with patience
we can see that all proper symmetric monomials can, after reindexing, be put in this form.
\end{proof}


\section{Steenrod algebra action}

We now  give a 
presentation of the cohomology of symmetric groups as algebras over
the mod-two Steenrod algebra $\A$.   

\begin{proposition}
The Steenrod squares $\sqi$ satisfy a Cartan formula with respect to transfer product.  That is
$\sqi (\alpha \tr \beta) = \sum_{j+k = i} \sq^{j} \alpha \tr \sq^{k} \beta$.
\end{proposition}

Because transfers are stable maps, they preserve Steenrod squares.  So this 
proposition is immediate from \refD{odot} and the external Cartan formula.
Because there are Cartan formulae with respect to both products, the cohomology
of $\coprod_{n} B \si_{n}$ is a Hopf ring over the Steenrod algebra.

The Steenrod algebra structure on all of the cohomology of $\coprod_{n} B \si_{n}$ is thus determined by 
the action on Hopf ring generators, and we consider the minimal set of $\gamma_{\ell, 2^{k}}$.
We introduce some notation to describe the action of Steenrod squares on these classes, using the
additive basis of gathered monomials for skyline diagrams from Section~\ref{S:comps}.

\begin{definition}
\begin{itemize}
\item The height of a gathered monomial is the largest of the algebraic degrees of its gathered blocks.  (The algebraic degree is the total number of Hopf ring generators cup-multiplied to give the gathered block.)
\item The effective scale of a gathered block, which is a product of $\gamma_{\ell,n}$'s, is the largest
$\ell$ which occurs.  The effective scale of a gathered monomial is the minimum of the effective 
scales of its constituent blocks.
\item We say a monomial is not full width if it is a non-trivial transfer product of some monomial with
some $1_{k}$.
\end{itemize}
\end{definition}

\begin{theorem}\label{T:steen}
$\sq^{i} \gamma_{\ell,2^{k}}$ is the sum of all full-width monomials of total degree $2^{k}(2^{\ell} - 1) + i$,
height one or two, and effective scale at least $\ell$.
\end{theorem}

We call such monomials the outgrowth monomials of  $\gamma_{\ell, 2^{k}}.$

For example, 
$$\sq^{3} \gamma_{2,4} = \gamma_{4,1} + \gamma_{3,1}\tr \gamma_{2,1} \gamma_{1,2} \tr \gamma_{2,1}
+ {\gamma_{2,1}}^{2} \tr \gamma_{2,1} \tr \gamma_{2,2}.$$

We can translate the conditions of \refT{steen} to our skyline diagrams, seeing that a Steenrod
square on $\gamma_{\ell,2^{k}}$ is represented by
the sum of all diagrams which are of full width, with at most two boxes stacked
on top of each other, and with the width of columns delineated by any of the vertical lines (of full height) at least $\ell$.
The example above translates to
$$
\begin{tabular}{ccccccc}
\raisebox{6pt}{\Large$\sq^{3}$(}\includegraphics[height=24pt]{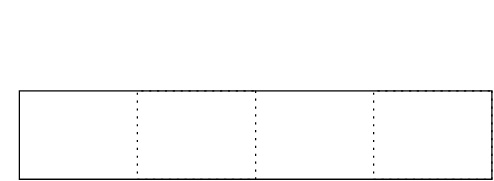} \raisebox{6pt}{\Large)}& \raisebox{10pt}{=} &\includegraphics[height=24pt]{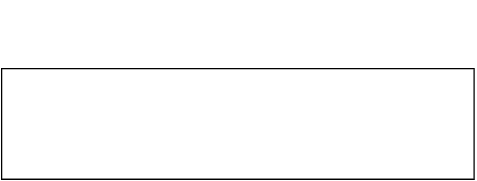} &\raisebox{10pt}{+}&  \includegraphics[height=24pt]{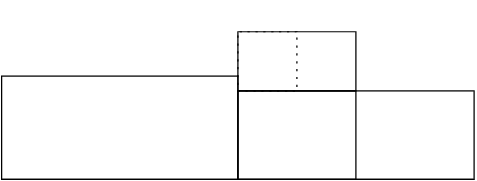}&\raisebox{10pt}{+}&\includegraphics[height=24pt]{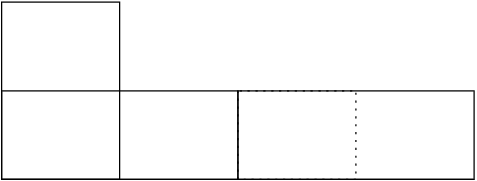}.\\
\end{tabular}
$$

We establish this theorem through restriction to suitable subgroups.
Recall that up to conjugation, the elementary abelian 
subgroups of $\si_{2^{n}}$ correspond to partitions of $2^{n}$ into powers
of two.  As before let $V_{n}$ denote the elementary abelian subgroup of $\si_{2^{n}}$ 
corresponding to the trivial partition.  
Include $\si_{2^{n - 1}} \times \si_{2^{n- 1}}$ in $\si_{2^{n}}$ in the standard way as in the definition
of the product in homology, so that the map on cohomology is the summand 
$\Delta_{2^{n-1}, 2^{n-1}}$ of the coproduct $\Delta$.  The following Lemma appears in \cite{MaMi79}.

\begin{lemma}\label{L:restrict}
 The sum of  restriction maps 
$$\hat{\rho} = {\rho}_{V_{n}} \oplus \Delta_{2^{n-1}, 2^{n-1}} : 
H^{*}(B\si_{2^{n}}) \to H^{*}(BV_{n}) \oplus H^{*}(B(\si_{2^{n- 1}} \times \si_{2^{n - 1}}))$$
is injective.
\end{lemma}

\begin{proof}
If we consider all elementary abelian subgroups of $B \si_{2^{n}}$ and thus all partitions of 
$2^{n}$ into powers of two, we see that other than the trivial partition such partitions must 
refine $2^{n} = 2^{n -1} + 2^{n  -1}$.   Thus the inclusions of
the corresponding 
elementary abelian subgroups factor up to conjugation through $\si_{2^{n -1}} \times \si_{2^{n-1}}$.    
The sum of restriction maps $\hat{\rho}$
is therefore injective because it factors the restriction to all elementary abelian subgroups.
\end{proof}

\begin{proof}[Proof of \refT{steen}]
We verify the equality of the theorem by verifying the agreement of  the restrictions of $\sqi \gamma_{\ell, 2^{k}}$ and the
sum of outgrowth monomials under $\rho_{V_{n}}$ and $\Delta_{2^{n-1}, 2^{n-1}}$ with $n = k + \ell$.   

\refC{maptodickson} states that the restriction of $\gamma_{\ell,2^{k}}$ to $V_{n}$ is just the Dickson class
$d_{k,\ell}$. In \cite{Hung91} Hu'ng calculated the Steenrod squares on Dickson classes as given by
$$\sqi d_{k, \ell}=  
\begin{cases}
d_{k', \ell'}  \;\; & i = 2^{k} - 2^{k'}\\
d_{k', \ell'} d_{k'', \ell''}  \;\;   & i = 2^{n} + 2^{k} - 2^{k'} - 2^{k''}, \;\; k' \leq k < k'' \\
{d_{k, \ell}}^{2} \;\; &i = 2^{k}(2^{\ell} - 1) \\
0 \;\;\;\;\;\;\;\; & {\text otherwise}.
\end{cases}
$$
By \refC{elem} the restriction to $V_{n}$ is zero for classes which are non-trivial transfer products and sends 
$\gamma_{\ell, 2^{k}}$ to $d_{k, \ell}$.  Thus 
the outgrowth monomials we need to consider in our formula for $\sq^{i} \gamma_{\ell, 2^{k}}$ are products
of one or two $\gamma_{\ell', 2^{k'}}$ with $\ell' + k' = n$ and (one) $\ell' > \ell$.  
Applying \refC{maptodickson} again we
see that this agrees with Hu'ng's calculations in Theorems A and B of \cite{Hung91}.

We show that the images of  $\sqi \gamma_{\ell, 2^{k}}$ and the
sum of outgrowth monomials under $\Delta_{2^{n-1}, 2^{n-1}}$ agree, by induction on $k$.  Since
we have already verified that all restrictions to $V_{n}$ agree then each inductive
step proves the theorem in that case.  If $k = 0$, the 
restriction of $\gamma_{\ell,1}$ is zero.  On the other hand, an outgrowth monomial of $\sqi \gamma_{\ell,1}$
is either zero  or a product $\gamma_{\ell, 1} \gamma_{\ell-k, 2^{k}}$ for $i = 2^{\ell} - 2^{k}$, which
restricts to zero because $\gamma_{\ell, 1}$ does.  

In general we have that
$$\Delta_{2^{n - 1}, 2^{n-1}} \gamma_{\ell, 2^{k}} = \gamma_{\ell,2^{k-1}} \otimes \gamma_{\ell, 2^{k-1}}.$$
We can thus apply the external Cartan formula to calculate that
$$\Delta_{2^{n - 1}, 2^{n-1}} \sqi \gamma_{\ell, 2^{k}} = 
\sum_{p+q = i} \sq^{p} \gamma_{\ell, 2^{k-1}} \otimes \sq^{q} \gamma_{\ell, 2^{k-1}},$$
which we understand by induction to be the sum of tensor products of two outgrowth monomials for 
$\gamma_{\ell, 2^{k-1}}$.   It thus suffices to show that the coproduct $\Delta_{2^{n-1}, 2^{n-1}}$
of the sum of all outgrowth monomials for $\gamma_{\ell, 2^{k}}$ is the sum of tensor products
of two outgrowth monomials for $\gamma_{\ell, 2^{k-1}}$.  This verification is straightforward using
\refP{gatheredcoprod}.    That such coproducts are given by sums of tensor products of
two such monomials is immediate since height and being full width are preserved 
by coproduct and effective scale can only increase on each factor.  All such tensor products
occur since we can form from $m \otimes m'$ a  monomial $M$ where if $\gamma_{\ell, p} \gamma_{\ell', p'}$
and $\gamma_{\ell, q} \gamma_{\ell, q'}$ are gathered blocks in $m$ and $m'$ respectively then 
$\gamma_{\ell, p+q} \gamma_{\ell', p' + q'}$ is a gathered block in $M$.  If $m$ and $m'$ are outgrowth monomials
for $\gamma_{\ell, 2^{k-1}}$ then $M$ will be for $\gamma_{\ell, 2^{k}}$ and $m \otimes m'$ will appear
in its coproduct.
\end{proof}

A geometric proof of this theorem might also be possible.  Because $\gamma_{\ell,2^{k}}$ are represented
by the varieties $\Gamma_{\ell, 2^{k}}$, the Wu formula implies that Steenrod operations on them are given by 
Stiefel-Whitney classes of the normal bundles to these varieties.  Since the $\Gamma_{\ell, 2^{k}}$ are defined by
equalities of coordinates, sections can be partially defined by perturbing those equalities and used for 
explicit computation.

In the previous section we revisited and extended the classical connection between the cohomology
of symmetric groups and Dickson algebras.  This connection persists when studying Steenrod
structure.   It has long been known that if filtered appropriately the cohomology of symmetric groups has as associated
graded a sum of exterior products of Dickson algebras, as algebras over the Steenrod algebra.  
From our point of view, we see this using the filtration by number of non-trivial transfer products
and then using the Cartan formulae we see that Steenrod squares do not increase filtration. 

\section{Cup product generators after Feshbach}\label{S:cupgens}

Feshbach gives in \cite{Fesh02}  a complete minimal 
set of ring generators for $H^{*}(B\si_{m}; \F_{2})$
along with relations which are minimal but not entirely explicit.
Using some results from the previous section, we can express
his generators in terms of our Hopf ring generators.  The 
combinatorics of even the generating set is somewhat involved.

\begin{definition}
A level-$n$ Dickson partition of $p$ is an equality 
$p = \sum_{k < n} t_{k} \left(2^{k}(2^{n-k} - 1)\right)$,
where at least one of the positive integers
$t_{k}$ is odd.  Consistent with \cite{Fesh02}, we denote such a partition
$\Lambda(n; {\bf t})$ or just $\Lambda$.
\end{definition}

\begin{definition}
Let $\Lambda$ be a level-$n$ Dickson partition.  Let
$2^{\ell}$ be the largest power of two which occurs twice in the dyadic expansion of the $t_{i}$'s
($\ell = -\infty$ if there is no such overlap), and  then define $\mu(\Lambda)$
to be $2^{\ell} + \sum 2^{d}$, where the sum is over all powers of 
two greater than $2^{\ell}$ which occur in
the dyadic expansion of some $t_{k}$.
In particular,  $\mu(\Lambda)$ is just the sum of the  $t_{k}$ if all powers of two in the dyadic 
expansion of the $t_{k}$ are distinct.
The maxwidth of $\Lambda$, denoted
$w(\Lambda)$, is defined as $2^{n} \mu(\Lambda)$.
\end{definition}

One of the main results of \cite{Fesh02} is the following.

\begin{theorem}
There is a minimal generating set of the cohomology of $B \si_{m}$ where the generators $v_{\Lambda}$ in dimension $p$ are indexed by Dickson partitions $\Lambda$ of $p$
with maxwidth less than $m$.
\end{theorem}

\begin{figure}
\begin{tabular}{| c | c | c|}     \hline
Level & Dickson partitions & Corresponding generators  (Hopf ring monomial and skyline names) \\ \hline
1     &  $1 = 1\cdot 1$ & $\gamma_{1,1} \tr 1_{10}$  \hspace{2cm} \includegraphics[width=50pt]{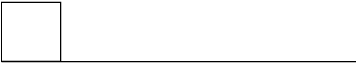} \\
 & $3 = 3\cdot 1$ &  $\gamma_{1,1}^{3} \tr 1_{10}$ \hspace{2cm} \includegraphics[width=50pt]{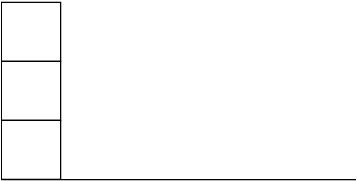} \\
&   $5 = 5\cdot 1$  &  $\gamma_{1,1}^{5} \tr 1_{10}$ \hspace{2cm} \includegraphics[width=50pt]{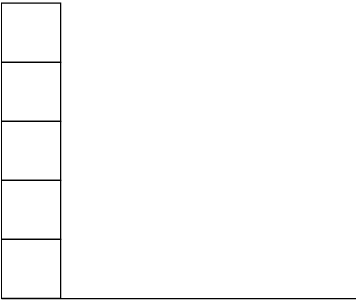}  \\  \hline
2      &  $2 = 1\cdot 2$ &  ${\gamma_{1,2}}  \tr 1_{8}$ \hspace{2cm}  \includegraphics[width=50pt]{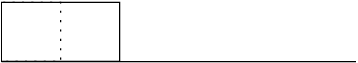} \\
&  $3 = 1 \cdot 3$ &  ${\gamma_{2,1}}  \tr 1_{8}$ \hspace{2cm}  \includegraphics[width=50pt]{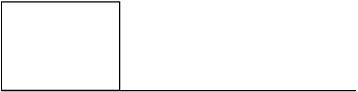} \\
& $5 = 1\cdot 2 + 1 \cdot 3$ & ${\gamma_{1,2}} {\gamma_{2,1}} \tr 1_{8}$ \hspace{1.5cm}  \includegraphics[width=50pt]{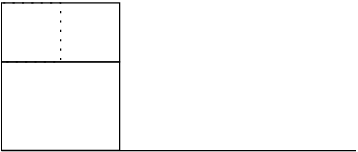} \\
& $6 = 3 \cdot 2 $ & ${\gamma_{1,2}}^{3}  \tr 1_{8}$ \hspace{2cm}  \includegraphics[width=50pt]{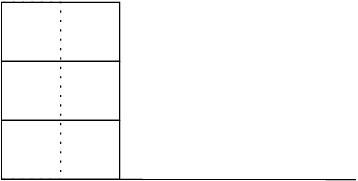} \\
&  $7 = 2 \cdot 2 + 1 \cdot 3$ & ${\gamma_{1,2}}^{2} {\gamma_{2,1}} \tr 1_{8}$ \hspace{1.5cm} \includegraphics[width=50pt]{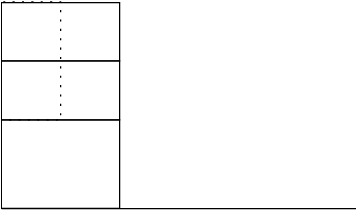} \\
& $8 = 1 \cdot 2 + 2 \cdot 3$ &  ${\gamma_{1,2}} {\gamma_{2,1}}^{2} \tr 1_{8}$ \hspace{1.5cm} \includegraphics[width=50pt]{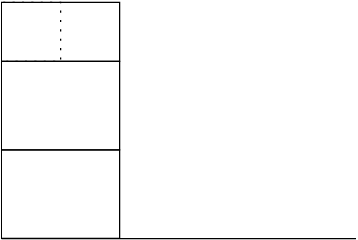}    \\ 
& $9 = 3 \cdot 3$ &  $ {\gamma_{2,1}}^{3} \tr 1_{8}$ \hspace{2cm}  \includegraphics[width=50pt]{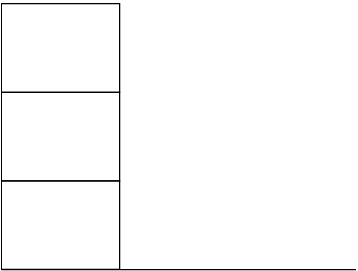}    \\  \hline
3      &    $4 = 1\cdot 4$ & ${\gamma_{1,4}}  \tr 1_{4}$ \hspace{2cm}  \includegraphics[width=50pt]{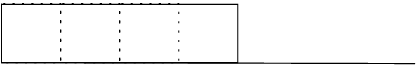} \\
& $6 = 1\cdot 6$ &   ${\gamma_{2,2}}  \tr 1_{4}$ \hspace{2cm} \includegraphics[width=50pt]{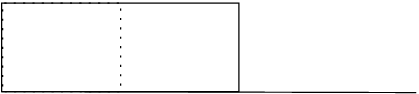} \\
& $7 = 1\cdot 7$ &  ${\gamma_{3,1}}  \tr 1_{4}$ \hspace{2cm} \includegraphics[width=50pt]{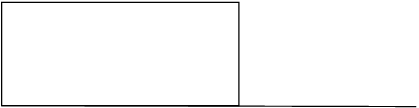}      \\ \hline 
\end{tabular}

\caption{A minimal generating set  under cup product for $H^{*}(B \si_{12})$.}
\end{figure}

These generators must be expressible in terms of our Hopf ring generators.

\begin{theorem}\label{T:LV}
The generator $v_{\Lambda}$ can be taken to be 
equal to $\prod_{k + \ell = n} {\gamma_{\ell, 2^k}}^{t_{k}} \tr
1_{m - 2^{k}}$.  Moreover, one obtains a generating set by replacing $\gamma_{\ell, 2^k}$'s
by any classes which pair non-trivially with $q_{k\cdot 0, \ell\cdot 1}$.
\end{theorem}

\begin{proof}[Proof of \refT{LV}]
First, we must recall Feshbach's characterization of multiplicative generators \cite{Fesh02}.
Recall that the elementary abelian $2$-subgroups
$W_{P}$ of $\si_{m}$ up to conjugacy 
are naturally indexed by partitions $P$ of $m$ as a sum of powers of two.  
If $\Lambda$
is a Dickson partition, we let $r_{\Lambda}$ be the smallest $r$ such that
$t_{r} \neq 0$.  We say that a partition $P$ of $m$ is $\mu$-subordinate to a level-$n$
Dickson partition $\Lambda$ if it contains a partition $\mu$ of the form 
$2^{n} = \sum_{j<r} s_{j} 2^{n-j}$.
The generator $v_{\Lambda}$ is characterized by its restriction to $W_{P}$ being
the sum $\sum_{\mu} v_{\mu, P}(\Lambda)$, where $\mu$ ranges over sub-partitions 
for which $P$ is $\mu$-subordinate to $\Lambda$ 
of $P$ and $v_{\mu, P}$ is the symmetrization of
$$1 \otimes 1 \otimes \cdots \otimes_{j} \left( \prod_{i}{d_{n-j, i-j}}^{t_{i}} \right)^{\otimes s_{j}} \otimes 1 \otimes \cdots \otimes 1,$$
where $d_{n-j, i-j}$ is the appropriate Dickson polynomial.

The class  $\prod_{k + \ell = n} {\gamma_{\ell, 2^k}}^{t_{k}} \tr 1_{m - 2^{n}}$ has 
these same restrictions.  The basic case of $m = 2^{n}$ is covered by \refC{maptodickson},
which says that the $\gamma_{\ell, 2^{k}}$ map to Dickson invariants $d_{k,\ell}$
in the cohomology of $V_{n}$.   The general case follows from the fact that
the inclusion $W_{P} = V_{n_{1}} \times \cdots \times V_{n_{q}}$ in $\si_{m}$
factors through the inclusion of $\si_{2^{n_{1}}} \times \cdots \times \si_{2^{n_{q}}}$,
which defines the product on homology.  We can thus use the coproduct
formula, namely \refL{AlgCoprod}, and the fact that the transfer coproduct is primitive to see that
$\prod_{k + \ell = n} {\gamma_{\ell, 2^k}}^{t_{k}} \tr 1_{m - 2^{k}}$ maps to  the cohomology 
of   $\si_{2^{n_{1}}} \times \cdots \times \si_{2^{n_{q}}}$ which in turn maps to the sum of $v_{\mu, P}$
as stated.
\end{proof}

We conjecture that alternate generating sets built from transfer products
alone (no cup products) of the $\gamma_{\ell,n}$, which might yield more tractable relations.


 \section{Stiefel-Whitney generators}\label{S:SW}

We now give an alternate presentation of this Hopf ring which uses Stiefel-Whitney
classes. Though it has a more complicated coproduct formula, such a presentation 
should be useful. The appearance of Stiefel-Whitney classes  forges another significant link between the
categories of finite sets and finite-dimensional vector spaces, showing that the mod-two cohomology
of automorphisms of finite sets is in a two-product sense generated by that of automorphisms
of vector spaces.

\begin{definition}
Let $w_{i,n} \in B\si_{n}$ be the pull-back of the $i$th  Stiefel-Whitney class
through the classifying map of the standard representation $\rho_{n}: \si_{n} \to O(n)$ given by
permutation matrices.
\end{definition}

\begin{remark}\label{SWrmk}
Stiefel-Whitney classes in symmetric groups have Poincar\'e dual representatives which are simple to describe in the configuration space model of $B \si_{n}$.  First replace $\uconf_{n}(\R^{\infty})$ by the homotopy equivalent subspace of configurations which are linearly independent.  The tautological bundle over $BO(n)$ pulls back to the bundle $E_{n}$ whose fiber over some configuration $\bf{x}$ is the vector space span $V_{{\bf{x}}}$ of the points $x_{i}$ in ${\bf{x}} = (x_{1}, \ldots, x_{n})/\sim$.   
Recall that $w_{i}$ of any bundle is the Poincar\'e dual to the locus where a generic collection of $n - i + 1$ sections becomes linearly dependent.  In this case, we may construct such sections by taking a standard basis $e_{1}, \cdots, e_{n-i+1}$ and projecting each one into $V_{{\bf{x}}}$.  Elementary linear algebra tells us that these projections will be dependent if and only if the projection of $V_{{\bf{x}}}$ onto $\R^{n-i+1} = {\rm Span}(e_{1}, \cdots, e_{n-i+1})$ is less than full rank.  

That is, the Poincar\'e dual of $w_{i}$ is the collection of unordered configurations  whose projection onto their first $n-i+1$ coordinates is not of full rank.  Thus for example $w_{1}$ records the linear dependence of a configuration of $n$ points in $\R^{\infty}$ when projected onto $\R^{n}$.  If we 
replace the bundle $E_{n}$ by $\bar{E_{n}}$, defined by taking the span of the vectors 
$x_{i} - x_{1}$, then $w_{2,4}$ is Poincar\'e dual to the locus of configurations of four points
in $\R^{\infty}$ which when projected onto $\R^{2}$ are collinear.    Thus, these Stiefel-Whitney classes are represented
by quadratic varieties, as opposed to the $\gamma_{\ell, n}$ which are represented by linear varieties.  But 
these Stiefel-Whitney varieties seem ``less singular.''

Through these Poincar\'e duals, we can explicitly see the pairings between Stiefel-Whitney
classes (and their cup and partition products) and polynomials in $q_{I}$ by counting 
intersections as in the figure above.  \end{remark}

\begin{figure}
$$\includegraphics[width=10cm]{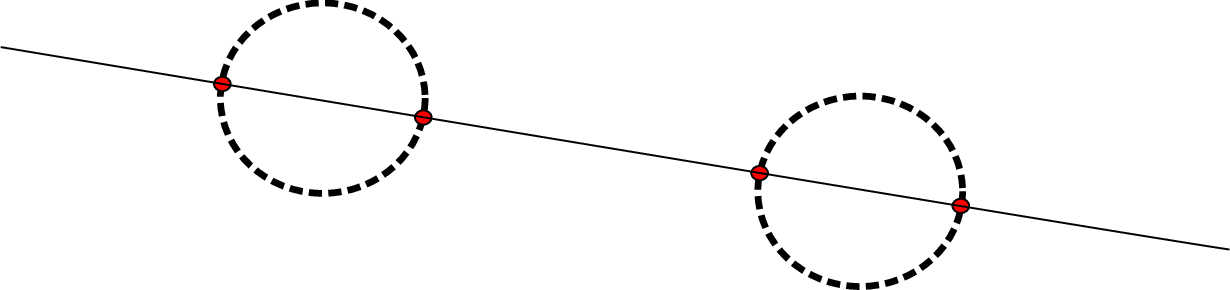}$$
\caption{An intersection  calculation, reflecting that $w_{2,4}$ pairs non-trivially with $q_{1} * q_{1}$.}
\end{figure}

\begin{definition}
Let $w(k, \ell) = w_{2^{k}(2^{\ell} - 1), 2^{k + \ell}}$.
\end{definition}

We will use \refC{GenCond} to show that the Stiefel-Whitney classes $w(k, \ell)$ generate
$H^{*}(\coprod_{n} B \si_{n})$ as a Hopf ring.  The needed calculation is a special case of what
is needed to understand the coproducts of these classes, so we treat the entire structure
at once.

\begin{proposition}\label{P:wiCalc}
Let ${\bf q}$ be a monomial in the $q_{I}$.  Then  $w_{i,n}$ evaluated on ${\bf q}$ 
is one if ${\bf q}$ is a product of classes  $q_{\ell \cdot 1}$ and $q_{0}$ and is zero otherwise.
\end{proposition}

\begin{proof}
 To set notation $H_{*}(\coprod_{n} BO(n))$ is the free polynomial ring on classes $b_{i}$ in
 degree $i$, with $i \geq 0$, which are in the homology of $BO(1) = \R P^{\infty}$.  Recall that 
  $w_{i,n} \in H^{i}(BO(n))$ evaluates to one on ${b_{0}}^{n-i}{b_{1}}^{i}$ and to zero on all other classes.  This calculation
 is easily established inductively using the fact that stably the coproduct of $w_{i}$ dual to the Pontrjagin product is $\sum w_{i-j} \otimes w_{j}$.  
 
 We use the fact that $\coprod_{n}B\rho_{n}$ is map of $E_{\infty}$-spaces (see for example II.7 of \cite{CLM76}) to see $q_{I} = q_{I}(\iota)$ maps to $q_{I}(b_{0})$.  Then we can use the structure of the homology of $\coprod_{n}BO(n)$ over the Kudo-Araki-Dyer-Lashof algebra, as computed by Kochman \cite{Koch73} and then Priddy \cite{Prid75}, determined by 
 \begin{equation}\label{E:KAonBO}
 q_{r}(b_{n}) = \sum_{i} \binom{r+i-1}{i} b_{n-i} b_{r+n+i}.   
 \end{equation}
Given that $w_{i,n}$ only pairs with ${b_{0}}^{n-i}{b_{1}}^{i}$ we call a monomial $m$, which 
is product of $b_{j}$'s, tainted if one of those $b_{j}$ has $j> 1$.  The basic observation is
that if $m$ is tainted then $q_{k}(m)$ is tainted for any $k$.  Indeed, using the Cartan formula
we see that $q_{k}(m)$ is sum of products of $q_{k_{i}}(b_{j_{i}})$.
But using \refE{KAonBO}, for the $j_{i} > 1$ the factor $q_{k_{i}}(b_{j_{i}})$ will be a sum of products of 
two $b$'s of total degree $2 j_{i} + k_{i} \geq 4$, so at least one must have degree greater than one.

If in $q_{i_{1}, \cdots i_{k}}$ we have $i_{k} > 1$, then because $q_{i_{k}}(b_{0}) = b_{0} b_{i_{k}}$
is tainted, so will be every term in $q_{i_{1}, \cdots i_{k}}(b_{0})$.  
Thus $w_{i,n}$ must evaluate trivially on any monomial which is a product of at least one such $q_{I}$. 

To see that  on the other hand $w_{i,n}$ does evaluate non-trivially on a monomial in the 
$b_{1,\ldots, 1}$, we calculate $q_{1, \ldots, 1}(b_{0})$.  We get that $q_{1}(b_{0}) = b_{0} b_{1}$,
and then that 
\begin{multline*}
q_{1,1}(b_{0}) = q_{1}(q_{1}(b_{0})) = q_{1}(b_{1} b_{1}) =
q_{1}(b_{0})  q_{0}(b_{1}) + q_{0}(b_{0}) q_{1}(b_{1})   
= b_{0} b_{1}^{3}  + b_{0}^{3} b_{3} + b_{0}^{2} b_{1} b_{2}.
\end{multline*}
In general, $q_{1, \ldots, 1}(b_{0})$ is equal 
 to $b_{0 }b_{1}^{2^{k}- 1}$ plus tainted monomials.
Because $B\rho_{*}$ is a map of rings, a product of such classes in degree $i$  will equal
${b_{0}}^{n-i}{b_{1}}^{i}$ plus tainted monomials, and thus
be evaluated non-trivially by $w_{i, n}$, completing the argument.
\end{proof}

We now develop the combinatorics
necessary to express the coproducts of Stiefel-Whitney classes.

\begin{definition}
A Dickson bi-partition of the pair $(k, \ell)$ is an equality of pairs of positive integers
$$\left( 2^{k}(2^{\ell} - 1), 2^{k+\ell} \right) =
\sum_{i} \left( 2^{k_{i}}(2^{\ell_{i}} - 1), 2^{k_{i} + \ell_{i}} \right).$$
Here we allow the trivial one-term partition, and we allow $k_{i}$ to be zero as well as $\ell_{i}$
to be zero when the corresponding $k_{i}$ is.  We manipulate such a partition as a set
$p = \{ (k_{i}, \ell_{i}) \}$, and sometimes emphasize the numbers being partitioned by
writing $p = p(k, \ell)$. 

 We say one bi-partition refines another if it is obtained by
substituting of some entry or entries by corresponding Dickson bi-partition(s).
\end{definition}

For example, because $(24, 32) = (4,8) + (6,8) + (14,16)$ we have the corresponding 
Dickson bi-parition $p(3,2) = \{ (2,1), (1,2), (1,3) \}$.  Because in turn of the equality 
$(4,8) = (0,2) + (1,2) + (3,4)$, we have that $q =  
\{ (1,0), (0,1), (1,1), (1,2), (1,3) \}$ refines $p$.

\begin{definition}
Let $\Pi_{k,\ell}$ denote the set consisting of Dickson bi-partitions $p$
expressed as an ordered union of two smaller partitions $p = p' \cup p''$,
each of which contains no repeated pairs of numbers.  Define a partial order
on $\Pi_{k,\ell}$ by $ p' \cup p''
\leq  q' \cup q''$ if $p'$ is a (possibly trivial) refinement of $q'$ and $p''$ of $q''$.

Let $\phi$ be the $\F_{2}$-valued function on $\Pi_{k,\ell}$ defined uniquely by
$$\sum_{p'\cup p''  \leq q' \cup q''} \phi(q' \cup q'')  = 1,$$ for any $p' \cup p'' \in \Pi_{k, \ell}$.
\end{definition}

In other words, the function $\phi$ is the inverse under convolution to the function 
which is one on all of $\Pi_{k, \ell}$.  Thus it could be determined by M\"obius inversion,
though we have not found that to be enlightening. 

\begin{theorem}\label{T:SWmain}
As a Hopf ring, $H^{*}(\coprod_{n} B \si_{n})$ is generated by Stiefel-Whitney classes
$w(k, \ell)$.

The transfer product is exterior, the antipode is the identity map, and there are no further relations other than those given
by Hopf ring distributivity and the fact that the product of classes on different components is zero.

The coproduct is given by
$$ \Delta w(k, \ell) = \sum_{p' \cup p'' \in \Pi_{k, \ell}}  \phi( p' \cup p'')
      \left( \bigodot_{(k_{i}, \ell_{i}) \in p'} w(k_{i}, \ell_{i}) \right) \bigotimes  
     \left( \bigodot_{(k_{j}, \ell_{j}) \in p''}  w(k_{j}, \ell_{j}) \right).$$
     
\end{theorem}

\begin{proof}
That these Stiefel-Whitney classes generate is now an immediate application of \refP{wiCalc}
to verify the hypothesis of \refC{GenCond}.  The lack of further relations and the additive basis 
follow from \refT{AlgMain} just as this \refC{GenCond} did.

The coproduct formula is verified by direct check using bialgebra structure.  
By \refP{wiCalc}, 
$w(k, \ell)$ evaluated on some non-trivial product $m * m'$ which is a monomial
will be non-zero if and only
if $m$ and $m'$ are products of ($q_{0}$'s and) $q_{1, \ldots, 1}$'s.  Such products are in one-to-one 
correspondence with the set $\Pi_{k, \ell}$.  Looking at only $m$, 
first we express each ${q_{1,...1}}^{n}$ uniquely
as a product of ${q_{1,...,1}}^{2^{k}} = q_{0, ........., 0, 1, ...., 1} = q(k, \ell)$, and then 
record the $(k, \ell)$ which appear.  For example,
$q_{0} {q_{1}}^{3} {q_{1,1}} {q_{1,1,1}}^{2} = q_{0} * q_{1} * q_{0,1} * q_{1,1} * q_{0,1,1,1}$ 
corresponds to $\{ (1,0),  (0,1), (1,1), (0,2), (1,3) \}$.  Call this
bijection $\beta$ from the set of monomials in $q_{1,...,1}$ to Dickson bi-partitions.

Applying \refP{wiCalc}, we find that not only does $m \otimes m'$ pair with 
$\bigodot_{(k_{i}, \ell_{i}) \in \beta(m)}w(k_{i}, \ell_{i}) \otimes \bigodot_{(k_{j}, \ell_{j}) \in \beta(m')}w(k_{j}, \ell_{j})$ but it also pairs with all similar products of Stiefel-Whitney
classes over $q' \cup q''$ which are refined by $\beta(m) \cup \beta(m')$.  Thus, if we take
the linear combination with coefficients given by $\phi$, that sum will pair
to one with $m \otimes m'$.
\end{proof}






\bibliographystyle{amsplain}
\bibliography{references}

\end{document}